\documentclass[11pt]{article}

\usepackage{amsmath}
\usepackage{amsfonts}
\usepackage{bm}
\usepackage{amssymb}
\usepackage{amsthm}
\usepackage{xcolor,url}
\usepackage{multirow}
\usepackage{graphicx}
\usepackage{epsfig}
\usepackage{subcaption}
\usepackage{epstopdf}
\usepackage{latexsym}
\usepackage{enumerate}
\usepackage{tabularx}
\usepackage{booktabs}
\usepackage{authblk}
\usepackage[left=1.3in, right=1.3in]{geometry}
\usepackage[bookmarksnumbered=true, citecolor=blue, colorlinks=true]{hyperref}
\usepackage[numbers,sort&compress]{natbib} 

\usepackage{mathtools}
\usepackage[nameinlink]{cleveref}
\crefname{equation}{}{} 

\usepackage{aliascnt}

\newtheorem{thm}{Theorem}[section]

\newcommand{\newaliasthm}[4]{%
  \newaliascnt{#1}{#3}%
  \newtheorem{#1}[#1]{#2}%
  \aliascntresetthe{#1}%
  \crefname{#1}{#2}{#4}%
  \Crefname{#1}{#2}{#4}%
}

\newaliasthm{cor}{Corollary}{thm}{Corollaries}
\newaliasthm{lem}{Lemma}{thm}{Lemmas}
\newaliasthm{prop}{Proposition}{thm}{Propositions}
\newaliasthm{defn}{Definition}{thm}{Definitions}
\newaliasthm{rem}{Remark}{thm}{Remarks}

\newaliasthm{example}{Example}{thm1}{Examples}

\numberwithin{equation}{section}

\newlength{\widebarargwidth}
\newlength{\widebarwidth}
\newlength{\widebarargheight}
\newlength{\widebarargdepth}

\DeclareMathOperator{\spn}{span}

\DeclareMathOperator{\filter}{filter}
\DeclareMathOperator{\extrap}{extrap}
\DeclareMathOperator{\hybrid}{hybrid}

\newcommand{\se}{\text{e}}
\newcommand{\si}{\text{i}}

\newcommand{\bR}{\mathbb{R}}

\newcommand{\bN}{\mathbb{N}}
\newcommand{\bZ}{\mathbb{Z}}

\newcommand{\bfH}{\mathbf{H}}

\newcommand{\vect}[1]{\bm{#1}}

\newcommand{\vc}{\vect{c}}

\newcommand{\veta}{\vect{\eta}}

\newcommand{\cG}{\mathcal{G}}

\newcommand{\cP}{\mathcal{P}}

\newcommand{\mtrx}[1]{\mathsf{#1}}

\newcommand{\mOmega}{\mtrx{\Omega}}

\providecommand{\keywords}[1]{\textbf{\text{Keywords.}} #1}

\begin{document}

\title{A hybrid reconstruction of piece-wise smooth functions from non-uniform Fourier data}

\author[1]{Guohui Song}
\author[2]{Congzhi Xia}
\affil[1]{Department of Mathematics and Statistics, Old Dominion University, Norfolk, VA 23529, USA. Email: gsong@odu.edu}
\affil[2]{Division of Science and Technology, Quincy University, Quincy, IL 62301, USA. Email: xiaco@quincy.edu}
\date{}

\maketitle

\begin{abstract}
  In this paper, we consider the problem of reconstructing piece-wise smooth functions from their non-uniform Fourier data. We first extend the filter method for uniform Fourier data to the non-uniform setting by using the techniques of admissible frames. We show that the proposed non-uniform filter method converges exponentially away from the jump discontinuities. However, the convergence rate is significantly slower near the jump discontinuities due to the Gibbs phenomenon. To overcome this issue, we combine the non-uniform filter method with a stable extrapolation method to recover the function values near the jump discontinuities. We show that the proposed hybrid method could achieve exponential accuracy uniformly on the entire domain. Numerical experiments are provided to demonstrate the performance of the proposed method.
\end{abstract}


\section{Introduction}
Fourier reconstruction is a fundamental problem in various applications such as signal processing, image processing, and medical imaging.
In many practical scenarios, the available Fourier data are often non-uniformly distributed due to physical constraints or limitations of the measurement devices.
The reconstruction of functions from non-uniform Fourier data has been extensively studied in the literature \cite{Benedetto2000,Christensen2003,Song2016,Adcock2014b,Gelb2012,Song2013b,Gelb2014a}.
However, most existing methods focus on the Fourier reconstruction of smooth functions, and the presence of jump discontinuities poses significant challenges due to the Gibbs phenomenon \cite{Gottlieb1997,Tadmor2007}.
We will consider the reconstruction of piece-wise smooth functions from non-uniform Fourier data in this paper.

When the Fourier data are uniformly distributed, various methods have been developed to mitigate the Gibbs phenomenon and achieve high accuracy in the reconstruction of piece-wise smooth functions.
These methods include spectral reprojection via the Gegenbauer polynomials \cite{Gottlieb1997,Gottlieb1992}, filter/mollifier methods \cite{Tadmor2002,Tadmor2005,Tanner2006}, generalized sampling techniques \cite{Adcock2014b}, and the algebraic system method \cite{Eckhoff1995,Eckhoff1998}.
In particular, the filter/mollifier method has been shown to achieve exponential accuracy when the parameters are chosen adaptively according to the distance from the reconstruction point to the nearest jump discontinuity \cite{Tadmor2002,Tadmor2005,Tanner2006}.
However, extending these methods to the non-uniform setting is not straightforward due to the lack of orthogonality of the non-uniform Fourier basis.

We will first extend the filter/mollifier method to the non-uniform Fourier reconstruction problem using the theory of frames \cite{Christensen2003}.
Specifically, we will employ the admissible frame approximation method \cite{Song2013b,Gelb2014a,Song2016} to handle the non-uniform Fourier data and then apply a filter in the physical domain to mitigate the Gibbs phenomenon.
We will show that this approach could achieve exponential accuracy away from the jump discontinuities.

However, the filter method may not perform well near the jump discontinuities due to the edge effects.
To address this issue, we will employ the stable extrapolation method \cite{Demanet2018} to recover the function values near the jump discontinuities.
The stable extrapolation method utilizes the analyticity of the function away from the jump discontinuities and approximates the function using Chebyshev polynomials.
It has been shown to achieve near-optimal accuracy in the extrapolation of analytic functions \cite{Demanet2018}.

Finally, we will combine the filter method and the stable extrapolation method to form a hybrid reconstruction method for piece-wise smooth functions from non-uniform Fourier data.
Specifically, we will first use the filter method to reconstruct the function values away from the jump discontinuities and then apply the stable extrapolation method to recover those values near the jump discontinuities based on the function values obtained from the filter method.
We show that the proposed hybrid method could achieve exponential accuracy \emph{uniformly} on the entire domain.

We point out that a hybrid reconstruction method combining the filter method and the Gegenbauer reconstruction method has been proposed in \cite{Gelb2000} for the uniform Fourier reconstruction of piece-wise smooth functions.
The proposed method in this paper could be viewed as an extension of the method in \cite{Gelb2000} to the non-uniform setting by extending the filter method to the non-uniform setting and replacing the Gegenbauer reconstruction method with the stable extrapolation method.
The non-uniform Fourier setting does not enjoy the nice correlation properties of the uniform Fourier basis with the Gegenbauer polynomials.
Additionally, the stable extrapolation method with Chebyshev polynomials is usually more stable than the Gegenbauer reconstruction method.
This makes the proposed hybrid method for non-uniform Fourier reconstruction more robust and easier to implement.

The rest of the paper is organized as follows.
In \Cref{sec:problem}, we describe the non-uniform Fourier reconstruction problem of piece-wise smooth functions and review the filter method for uniform Fourier reconstruction.
In \Cref{sec:filter}, we extend the filter method to the non-uniform Fourier reconstruction problem using the theory of frames and prove that it could achieve exponential accuracy away from the jump discontinuities.
In \Cref{sec:hybrid}, we introduce the stable extrapolation method and combine it with the filter method to form a hybrid reconstruction method for piece-wise smooth functions from non-uniform Fourier data.
We also present the hybrid method that could achieve exponential accuracy uniformly on the entire domain.
In \Cref{sec:numerical}, we provide numerical experiments to demonstrate the performance of the proposed hybrid method.
Finally, we conclude the paper in \Cref{sec:conclusion}.

\section{Non-uniform Fourier reconstruction of piece-wise smooth functions}\label{sec:problem}
We will first describe the non-uniform Fourier reconstruction problem of piece-wise smooth functions. Suppose we have a piece-wise smooth function \(f: \bR \rightarrow \bR\) that is supported on \([0,1]\). For simplicity of presentation, we assume that \(f\) has a single jump discontinuity \(\xi \in (0,1)\) and is infinitely smooth on \((-\infty, \xi)\) and \((\xi, \infty)\) respectively. Our proposed hybrid reconstruction method could be easily extended to the case when \(f\) has multiple jump discontinuities. We will assume the jump discontinuity \(\xi\) is known a priori. Otherwise, we could use existing edge detection methods \cite{Cochran2013,Gelb2016,Gelb2011,Stefan2012,Xiao2023} to estimate the location of the jump discontinuity from the given non-uniform Fourier data.
We are given finitely many non-uniform Fourier samples of \(f\) in the following form:
\begin{align*}
  \hat{f}(\lambda_j) = \langle f, \psi_j \rangle = \int_0^1 f(x) \psi_j^*(x) dx, \quad j = -m, \cdots, m,
\end{align*}
where \(\psi_j(x) = \se^{2\pi i \lambda_j x}\) for some non-uniformly distributed \(\{\lambda_j\}_{j=-m}^m\) in \(\bR\) and \(\psi_{j}^{*}\) is the complex conjugate of \(\psi_{j}\). We would like to reconstruct \(f\) from the above non-uniform Fourier data.

A key challenge in the Fourier reconstruction problem is due to the presence of the jump discontinuity, which causes the well-known Gibbs phenomenon in the standard Fourier partial sum reconstruction.
Additionally, the non-uniform distribution of the Fourier samples further complicates the reconstruction problem.
We will first extend the filter method for uniform Fourier reconstruction to the non-uniform setting and then introduce a stable extrapolation method to recover the function values near the jump discontinuity.
Finally, we will combine these two methods to form a hybrid reconstruction method for piece-wise smooth functions from non-uniform Fourier data.

We first review the filter method \cite{Tadmor2002,Tadmor2005,Tadmor2007} for uniform Fourier reconstruction. When given the uniform Fourier coefficients \(\{\hat{f}(j)\}_{|j|\le m}\), the filter method reconstructs \(f\) by
\begin{align*}
  Q_m^\sigma f(x) = \sum_{|j|\le m} \sigma\biggl(\frac{|j|}{m}\biggr) \hat{f}(j) \se^{2\pi i j x},
\end{align*}
where \(\sigma\) is a filter function that satisfies certain admissible conditions. There are various types of filters investigated in the literature. The classical choices of filters with compact support could achieve root exponential accuracy \cite{Tadmor2002,Tadmor2005}. A non-compact supported filter with a joint time-frequency localization utilizing the Hermite Distributed Approximating Functionals (HDAFs) has been shown in \cite{Tanner2006} to obtain exponential accuracy. We will consider the HDAF-based filter in this paper. Specifically, the \(p\)-th order HDAF is defined as
\begin{align*}
  \bfH_{p,\gamma}(x) = \frac{1}{\gamma} e^{-(\frac{x}{\sqrt{2}\gamma})^2} \sum_{l=0}^p \frac{(-1)^l}{4^l l!} H_{2l}\biggl(\frac{x}{\sqrt{2}\gamma}\biggr),
\end{align*}
where \(\gamma\) is a free localization parameter and \(H_{2l}\) is the Hermite polynomial of order \(2l\). The corresponding filter function is defined by the Fourier transform of the HDAF:
\begin{align*}
  \sigma_{p,\gamma}(w) = \frac{1}{\sqrt{2}} \int_{-\infty }^{\infty }\bfH_{p,\gamma}(x)  e^{-\pi i w x} dx = \se^{-\frac{(w\gamma)^{2}}{2}}\sum_{l=0}^p \frac{1}{l! 2^l} (w \gamma)^{2l}, \quad w \in \bR.
\end{align*}

The filter method could also be viewed in an equivalent way of applying a mollifier in the physical domain. In particular, for the above HDAF-based filter \(\sigma_{p,\gamma}\), we have \cite{Tanner2006}
\begin{align*}
  Q_m^{\sigma_{p,\gamma}} f(x) = (\rho_{p,\gamma} * Q_m f)(x), \quad x\in \bR,
\end{align*}
where the mollifier \(\rho_{p,\gamma}\) is defined by
\begin{align}\label{eq:mollifier}
  \rho_{p,\gamma}(x) = \sum_{j=-\infty}^\infty \bfH_{p,\gamma}(m(x+2j)), \quad x\in \bR,
\end{align}
and \(Q_m f\) is the orthogonal projection of \(f\) onto the space \(\cG_m = \spn\{\phi_j: |j|\le m\}\) with \(\phi_j(x) = \se^{2\pi i j x}\) being the classical Fourier basis.

The choices of the parameters \(p\) and \(\gamma\) are crucial to the performance of the filter method.
It was shown in \cite{Tanner2006} that if they are chosen adaptively according to the distance from \(x\) to the nearest jump discontinuity, then the filter method could achieve exponential accuracy.

We will try to extend the above filter method to the non-uniform Fourier reconstruction problem.
However, it is not straightforward to do so since the non-uniform Fourier samples at \(\{\lambda_j\}_{j=-m}^m\) do not form an orthogonal basis for \(L^2[0,1]\).
We will use the theory of frames to handle the non-uniform Fourier samples.
Specifically, we will assume that the non-uniform Fourier samples \(\{\lambda_j\}_{j\in \bZ}\) forms a Fourier frame for \(L^2[0,1]\).
That is, there exist two positive constants \(A\) and \(B\) such that
\begin{align*}
  A \|f\|_2 \leq \sum_{j\in \bZ} |\hat{f}(\lambda_j)|^2 \leq B \|f\|_2
\end{align*}
for all \(f \in L^2[0,1]\), where \(\|\cdot\|_2\) denotes the \(L^2\) norm. We will present the filter method for non-uniform Fourier reconstruction in the next section.

\section{Filter method for non-uniform Fourier reconstruction}\label{sec:filter}
We will develop the filter method for non-uniform Fourier reconstruction in this section. We point out that the \emph{admissible frame approximation} method has been introduced in \cite{Song2013b,Gelb2014a,Song2016} to reconstruct a smooth function from its non-uniform Fourier data. It has been further employed in \cite{Gelb2016} to detect the jump discontinuities of a piece-wise smooth function. However, the admissible frame approximation method is not directly applicable to the reconstruction of piece-wise smooth functions due to the Gibbs phenomenon. We will combine the filter method and the admissible frame approximation method to reconstruct piece-wise smooth functions from non-uniform Fourier data.

We give a brief review of the admissible frame approximation method \cite{Song2013b,Gelb2014a,Song2016} in the following. Suppose \(\{\psi_{j}(x) = \se^{2 \pi i \lambda_j x}\}_{j\in \bZ}\) is a frame for \(L^2[0,1]\). Then, for any smooth function \(f\in L^2[0,1]\), we have the following approximation based on its non-uniform Fourier samples \(\{\hat{f}(\lambda_j)\}_{j=-m}^m\):
\begin{align*}
  W_{n}^{-1} Q_n S_m f,
\end{align*}
where \(S_m f = \sum_{|j|\le m} \hat{f}(\lambda_j) \psi_j\), \(Q_n\) is the orthogonal projection onto \(\cG_n = \spn\{\phi_j: |j|\le n\}\) with \(\phi_j(x) = \se^{2\pi i j x}\) being the classical Fourier basis, and \(W_n\) is the restriction of the frame operator \(W = Q_n S_m|_{\cG_n}\) on \(\cG_n\). It was shown in \cite{Song2013b} that if the frame \(\{\psi_j\}_{j\in \bZ}\) satisfies certain localization conditions and the parameter \(n\) is chosen appropriately, then the above approximation could achieve high order convergence for smooth functions. In particular, if the Fourier frame is \emph{admissible} with respect to the Fourier basis \(\{\phi_j\}_{j\in \bZ}\) in the sense that there exist positive constant \(c_0\) such that
\begin{align}\label{eq:admissible}
  |\langle \psi_j, \phi_l \rangle| \leq c_0 (1+|j-l|)^{-1}, \quad j,l \in \bZ,
\end{align}
and we choose \(n\) as
\begin{align}\label{eq:admissible_n}
  n = \frac{A m}{A + 2 c_{0}^{2}},
\end{align}
with \(A\) being the lower frame bound, then we have the following results \cite{Song2013b}: for any \(f\in L^2[0,1]\),
\begin{align}\label{eq:admissible_error}
  \|Q_{n}f - W_n^{-1} Q_n S_m f\|_2 \leq \frac{2B}{A} \left\lVert f - Q_{n}f\right\rVert_2.
\end{align}
The above results indicate that the convergence rate of the admissible frame approximation is determined by the convergence rate of the projection \(Q_n\). In particular, if \(f\) is smooth, then we have the spectral convergence. However, if \(f\) is piece-wise smooth, then the convergence rate of \(Q_n\) is only first order due to the Gibbs phenomenon. We will combine the filter method and the admissible frame approximation method to reconstruct piece-wise smooth functions from non-uniform Fourier data in the following.

We will employ the admissible frame approximation method to construct a filtered version of the function. Specifically, for any \(x\in [0,1]\), we will approximate \(f(x)\) by
\begin{align}\label{eq:filter_reconstruction}
  f_{m}^{\filter}(x) = W_{n}^{-1} Q_n S_m (\rho_{p_{x}, \gamma_{x}} *  f)(x),
\end{align}
where \(\rho_{p_{x}, \gamma_{x}}\) is the mollifier defined in \cref{eq:mollifier} with parameters \(p_{x}\) and \(\gamma_{x}\) chosen adaptively according to the distance from \(x\) to the jump discontinuity \(\xi\). We will discuss the choice of \(p_{x}\) and \(\gamma_{x}\) later. Note that we could compute the function \(f_{m}^{\filter}\) from the given nonuniform Fourier data \(\{\hat{f}(\lambda_{j})\}_{j=-m}^{m}\). Specifically, as shown in \cite{Song2013b}, we have
\begin{align*}
  f_{m}^{\filter}(x) = \sum_{l=-n}^{n}c_{x,l}\phi_{l}(x),
\end{align*}
where the coefficients vector \(\vc_{x}=\left(c_{x,l}:-n\leq l\leq n\right)\) is given by \(\vc_{x} = \mOmega^{\dagger}\veta_{x}\), \(\mOmega^{\dagger}\) is the Moore-Penrose pseudo-inverse of the matrix \(\mOmega = \left[\langle \psi_j, \phi_l \rangle: -m\leq j\leq m, -n\leq l \leq n\right]\), and \(\veta_{x} = \left(\hat{\rho}_{p_{x}, \gamma_{x}}(\lambda_{j}) \hat{f}(\lambda_{j}): -m\leq j\leq m\right) \).

We next show that the above filter reconstruction could achieve exponential accuracy away from the jump discontinuity. In particular, for any \(x\in [0,1]\), we will analyze the error \(E(x) = |f(x) - f_{m}^{\filter}(x)|\) through decomposing it into the following two parts:
\begin{align*}
  E_{1}(x)  = |f(x) - (\rho_{p_{x}, \gamma_{x}} * f)(x)|
\end{align*}
and
\begin{align*}
  E_{2}(x)  = |(\rho_{p_{x}, \gamma_{x}} * f)(x) - f_{m}^{\filter}(x)|.
\end{align*}
The first part \(E_{1}(x)\) has been analyzed in \cite{Tanner2006} and it could achieve exponential accuracy when the parameters \(p_{x}\) and \(\gamma_{x}\) are chosen appropriately. Indeed, for any \(x \in [0,1]\), let \(d(x) = |x - \xi|\) be the distance from \(x\) to the jump discontinuity \(\xi\). If we choose
\begin{align}\label{eq:gammax_px}
  \gamma_{x} = \sqrt{\alpha d(x)m}, \quad p_x = \kappa d(x) m,
\end{align}
where \(\alpha\) and \(\kappa\) are some positive constants such that \(\alpha \kappa < \frac{1}{2\ln(1+\sqrt{2})}\), then we have the following exponential accuracy \cite{Tanner2006}:
\begin{align}\label{eq:E1}
  E_{1}(x) \le C_f m^{9/4} \tau^{-m d(x)},
\end{align}
where \(\tau = \min(\eta_f^\kappa, \frac{\se^{1/4\alpha}}{(1+\sqrt{2})^\kappa})\), \(\eta_f > 1\) is a constant depending on the analytic extension of \(f\) to the complex plane, and \(C_f\) is a constant depending on \(f\).

It remains to analyze the second part \(E_{2}(x)\). To this end, we first estimate \(\rho_{p_{x}, \gamma_{x}} * f - Q_{n}(\rho_{p_{x}, \gamma_{x}} * f)\), which will play a key role in the analysis of \(E_{2}(x)\). We have the following result.
\begin{lem}\label{lem:mollifier_projection}
  For any \(x\in [0,1]\), if \(\frac{n^{2}\gamma_{x}^{2}}{2m^{2}} \geq p_x\), then we have
  \begin{align*}
    \| \rho_{p_{x}, \gamma_{x}} * f - Q_{n}(\rho_{p_{x}, \gamma_{x}} * f)\|_{\infty} \leq 2 \left\lVert f\right\rVert_{\infty} \frac{n}{p_{x}!} \se^{-\frac{n^{2}\gamma_{x}^{2}}{2m^{2}}} \left(\frac{n\gamma_{x}}{\sqrt{2}m}\right)^{2p_{x}}
  \end{align*}
  and
  \begin{align*}
    \| \rho_{p_{x}, \gamma_{x}} * f - Q_{n}(\rho_{p_{x}, \gamma_{x}} * f)\|_{2} \leq  \left\lVert f\right\rVert_{\infty} \frac{\sqrt{2n}}{p_{x}!} \se^{-\frac{n^{2}\gamma_{x}^{2}}{2m^{2}}} \left(\frac{n\gamma_{x}}{\sqrt{2}m}\right)^{2p_{x}}.
  \end{align*}
\end{lem}
\begin{proof}
  Note that
  \begin{align*}
    \rho_{p_{x}, \gamma_{x}} * f - Q_{n}(\rho_{p_{x}, \gamma_{x}} * f) = \sum_{|j| >n} \langle \rho_{p_{x}, \gamma_{x}} * f, \phi_j \rangle \phi_j.
  \end{align*}
  It has been shown in \cite{Tanner2006} that for any \(j\in \bZ\),
  \begin{align*}
    |\langle \rho_{p_{x}, \gamma_{x}} * f, \phi_j \rangle| \leq \left\lVert f\right\rVert_{\infty} \se^{-\frac{j^{2}\gamma_{x}^{2}}{2m^{2}}}\sum_{l=0}^{p_{x}} \frac{1}{l!} \left(\frac{j^{2}\gamma_{x}^{2}}{2m^{2}}\right)^{l}.
  \end{align*}
  When \(j\geq n\), we have \(\frac{j^{2}\gamma_{x}^{2}}{2m^{2}}\geq p_{x}\) and the last term in the above summation dominates the others. Thus, we have for any \(j\geq n\),
  \begin{align}\label{eq:mollifier_Fourier_coeff}
    |\langle \rho_{p_{x}, \gamma_{x}} * f, \phi_j \rangle| \leq \left\lVert f\right\rVert_{\infty} \se^{-\frac{j^{2}\gamma_{x}^{2}}{2m^{2}}} \frac{1}{p_{x}!} \left(\frac{j^{2}\gamma_{x}^{2}}{2m^{2}}\right)^{p_{x}}.
  \end{align}

  We first estimate the \(L^{\infty}\) norm. For any \(n\geq 1\), we have
  \begin{align*}
    \| \rho_{p_{x}, \gamma_{x}} * f - Q_{n}(\rho_{p_{x}, \gamma_{x}} * f)\|_{\infty} & \leq \sum_{|j|>n} |\langle \rho_{p_{x}, \gamma_{x}} * f, \phi_j \rangle|.
  \end{align*}
  It follows that
  \begin{align*}
    \| \rho_{p_{x}, \gamma_{x}} * f - Q_{n}(\rho_{p_{x}, \gamma_{x}} * f)\|_{\infty} & \leq 2 \left\lVert f\right\rVert_{\infty} \sum_{j>n} \se^{-\frac{j^{2}\gamma_{x}^{2}}{2m^{2}}} \frac{1}{p_{x}!} \left(\frac{j^{2}\gamma_{x}^{2}}{2m^{2}}\right)^{p_{x}}.
  \end{align*}
  When \(\frac{j^{2}\gamma_{x}^{2}}{2m^{2}}\geq p_{x}\), the integrand in the above summation is decreasing with respect to \(j\). Thus, we have
  \begin{align*}
    \| \rho_{p_{x}, \gamma_{x}} * f - Q_{n}(\rho_{p_{x}, \gamma_{x}} * f)\|_{\infty} & \leq 2 \left\lVert f\right\rVert_{\infty} \int_{n}^{\infty} \se^{-\frac{t^{2}\gamma_{x}^{2}}{2m^{2}}} \frac{1}{p_{x}!} \left(\frac{t^{2}\gamma_{x}^{2}}{2m^{2}}\right)^{p_{x}} dt.
  \end{align*}
  By a change of variable \(u = \frac{t\gamma_{x}}{\sqrt{2}m}\), we have
  \begin{align}\label{eq:mollifier_projection_Linf}
    \| \rho_{p_{x}, \gamma_{x}} * f - Q_{n}(\rho_{p_{x}, \gamma_{x}} * f)\|_{\infty} & \leq 2 \left\lVert f\right\rVert_{\infty} \frac{\sqrt{2}m}{\gamma_{x} p_{x}!} \int_{\frac{n\gamma_{x}}{\sqrt{2}m}}^{\infty} \se^{-u^{2}} u^{2p_{x}} du.
  \end{align}
  Applying integration by parts on the last integral, we obtain
  \begin{align*}
    \int_{\frac{n\gamma_{x}}{\sqrt{2}m}}^{\infty} \se^{-u^{2}} u^{2p_{x}} du = \frac{1}{2} \se^{-\frac{n^{2}\gamma_{x}^{2}}{2m^{2}}} \left(\frac{n\gamma_{x}}{\sqrt{2}m}\right)^{2p_{x}-1} + \frac{2p_{x}-1}{2} \int_{\frac{n\gamma_{x}}{\sqrt{2}m}}^{\infty} \se^{-u^{2}} u^{2p_{x}-2} du.
  \end{align*}
  Since \(\frac{n\gamma_{x}}{\sqrt{2}m}\geq  1\), we have
  \begin{align*}
    \int_{\frac{n\gamma_{x}}{\sqrt{2}m}}^{\infty} \se^{-u^{2}} u^{2p_{x}} du \leq  \frac{1}{2} \se^{-\frac{n^{2}\gamma_{x}^{2}}{2m^{2}}} \left(\frac{n\gamma_{x}}{\sqrt{2}m}\right)^{2p_{x}-1} + \frac{2p_{x}-1}{2}\left(\frac{n^{2}\gamma_{x}^{2}}{2m^{2}}\right)^{-1} \int_{\frac{n\gamma_{x}}{\sqrt{2}m}}^{\infty} \se^{-u^{2}} u^{2p_{x}-2} du.
  \end{align*}
  When \(\frac{n^{2}\gamma_{x}^{2}}{2m^{2}}\geq p_{x}\), we have
  \begin{align*}
    \int_{\frac{n\gamma_{x}}{\sqrt{2}m}}^{\infty} \se^{-u^{2}} u^{2p_{x}} du \leq \frac{\frac{n^{2}\gamma_{x}^{2}}{2m^{2}}}{2\frac{n^{2}\gamma_{x}^{2}}{2m^{2}} - 2p_{x} + 1} \se^{-\frac{n^{2}\gamma_{x}^{2}}{2m^{2}}} \left(\frac{n\gamma_{x}}{\sqrt{2}m}\right)^{2p_{x}-1} \leq \se^{-\frac{n^{2}\gamma_{x}^{2}}{2m^{2}}} \left(\frac{n\gamma_{x}}{\sqrt{2}m}\right)^{2p_{x}+1}.
  \end{align*}
  Substituting the above estimate into \eqref{eq:mollifier_projection_Linf}, we obtain
  \begin{align*}
    \| \rho_{p_{x}, \gamma_{x}} * f - Q_{n}(\rho_{p_{x}, \gamma_{x}} * f)\|_{\infty} \leq 2 \left\lVert f\right\rVert_{\infty} \frac{n}{p_{x}!} \se^{-\frac{n^{2}\gamma_{x}^{2}}{2m^{2}}} \left(\frac{n\gamma_{x}}{\sqrt{2}m}\right)^{2p_{x}}.
  \end{align*}

  We next estimate the \(L^{2}\) norm. For any \(n\geq 1\), we have
  \begin{align*}
    \| \rho_{p_{x}, \gamma_{x}} * f - Q_{n}(\rho_{p_{x}, \gamma_{x}} * f)\|_{2}^2 & = \sum_{|j| >n} |\langle \rho_{p_{x}, \gamma_{x}} * f, \phi_j \rangle|^2.
  \end{align*}
  Substituting \eqref{eq:mollifier_Fourier_coeff} into the above equation, we obtain
  \begin{align*}
    \| \rho_{p_{x}, \gamma_{x}} * f - Q_{n}(\rho_{p_{x}, \gamma_{x}} * f)\|_{2}^2 & \leq 2 \left\lVert f\right\rVert_{\infty}^2 \sum_{j>n} \se^{-\frac{j^{2}\gamma_{x}^{2}}{m^{2}}} \frac{1}{(p_{x}!)^2} \left(\frac{j^{2}\gamma_{x}^{2}}{2m^{2}}\right)^{2p_{x}}.
  \end{align*}
  By a similar argument as the estimate for the \(L^{\infty}\) norm, we have
  \begin{align*}
    \| \rho_{p_{x}, \gamma_{x}} * f - Q_{n}(\rho_{p_{x}, \gamma_{x}} * f)\|_{2}^2 \leq 2 \left\lVert f\right\rVert_{\infty}^{2} \frac{n}{(p_{x}!)^2} \se^{-\frac{n^{2}\gamma_{x}^{2}}{m^{2}}} \left(\frac{n\gamma_{x}}{\sqrt{2}m}\right)^{4p_{x}},
  \end{align*}
  which implies the desired estimate for the \(L^{2}\) norm.
\end{proof}

We next estimate the error \(E_{2}(x)\). We have the following result.
\begin{prop}\label{prop:E2}
  Suppose the Fourier frame \(\{\psi_j\}_{j\in \bZ}\) is admissible with respect to the Fourier basis \(\{\phi_j\}_{j\in \bZ}\) as defined in \cref{eq:admissible}. For any \(x\in [0,1]\), if we choose \(n\) as \cref{eq:admissible_n}, \(p_x\) and \(\gamma_x\) such that \(\frac{n^{2}\gamma_{x}^{2}}{2m^{2}} \geq p_x\), then we have
  \begin{align*}
    E_{2}(x) \leq \frac{A+2B}{A} \frac{2n+1}{p_{x}!} \se^{-\frac{n^{2}\gamma_{x}^{2}}{2m^{2}}} \left(\frac{n\gamma_{x}}{\sqrt{2}m}\right)^{2p_{x}}.
  \end{align*}
\end{prop}
\begin{proof}
  We first decompose \(E_{2}(x)\) as follows:
  \begin{align*}
    E_{2}(x)  \leq  |(\rho_{p_{x}, \gamma_{x}} * f)(x) - Q_{n}(\rho_{p_{x}, \gamma_{x}} * f)(x)| + |Q_{n}(\rho_{p_{x}, \gamma_{x}} * f)(x) - f_{m}^{\filter}(x)|.
  \end{align*}
  The first term has been estimated in \Cref{lem:mollifier_projection}. For the second term, we note that both \(Q_{n}(\rho_{p_{x}, \gamma_{x}} * f)\) and \(f_{m}^{\filter}\) belong to the space \(\cG_n\). Thus, we have
  \begin{align*}
    |Q_{n}(\rho_{p_{x}, \gamma_{x}} * f)(x) - f_{m}^{\filter}(x)| & \leq \sqrt{2n+1} \|Q_{n}(\rho_{p_{x}, \gamma_{x}} * f) - f_{m}^{\filter}\|_2.
  \end{align*}
  By the definition of \(f_{m}^{\filter}\) \cref{eq:filter_reconstruction} and the error estimate \cref{eq:admissible_error}, we have
  \begin{align*}
    \|Q_{n}(\rho_{p_{x}, \gamma_{x}} * f) - f_{m}^{\filter}\|_2 \leq \frac{2B}{A} \|\rho_{p_{x}, \gamma_{x}} * f - Q_{n}(\rho_{p_{x}, \gamma_{x}} * f)\|_2,
  \end{align*}
  where \(A\) and \(B\) are the lower and upper frame bounds of the Fourier frame \(\{\psi_j\}_{j\in \bZ}\). It then follows that
  \begin{align*}
    |Q_{n}(\rho_{p_{x}, \gamma_{x}} * f)(x) - f_{m}^{\filter}(x)| & \leq \sqrt{2n+1} \frac{2B}{A} \|\rho_{p_{x}, \gamma_{x}} * f - Q_{n}(\rho_{p_{x}, \gamma_{x}} * f)\|_2.
  \end{align*}
  It implies
  \begin{align*}
    E_{2}(x)  \leq & \|(\rho_{p_{x}, \gamma_{x}} * f)- Q_{n}(\rho_{p_{x}, \gamma_{x}} * f)\|_{\infty}+ \sqrt{2n+1} \frac{2B}{A} \|\rho_{p_{x}, \gamma_{x}} * f - Q_{n}(\rho_{p_{x}, \gamma_{x}} * f)\|_2.
  \end{align*}
  The desired result then follows from applying \Cref{lem:mollifier_projection} to the above inequality.
\end{proof}

We are now ready to combine the estimates of \(E_{1}(x)\) and \(E_{2}(x)\) to obtain the convergence result of the filter reconstruction \(f_{m}^{\filter}\).
\begin{thm}\label{thm:filter_reconstruction}
  Suppose the Fourier frame \(\{\psi_j\}_{j\in \bZ}\) is admissible with respect to the Fourier basis \(\{\phi_j\}_{j\in \bZ}\) as defined in \cref{eq:admissible}. For any \(x\in [0,1]\), if we choose \(n\) as in \cref{eq:admissible_n}, \(p_x\) and \(\gamma_x\) as in \cref{eq:gammax_px}, where the positive constants \(\alpha\) and \(\kappa\) there satisfy
  \begin{align*}
    \alpha \kappa < \frac{1}{2\ln(1+\sqrt{2})} \quad \text{and} \quad \frac{\kappa}{\alpha} \leq \frac{A^{2}}{2(A + 2 c_{0}^{2})^{2}},
  \end{align*}
  then we have
  \begin{align*}
    |f(x) - f_{m}^{\filter}(x)| \leq C_f m^{9/4} \tau^{-m d(x)} + \frac{(A+2B)}{A \sqrt{2\pi \kappa d(x) m}} \left(\frac{2Am}{A + 2 C_{0}^{2}}+1\right) \se^{-(\theta - 1 - \ln(\theta))\kappa d(x) m}.
  \end{align*}
  where \(d(x) = |x-\xi|\) is the distance from \(x\) to the jump discontinuity \(\xi\), \(\theta= \frac{\alpha A^{2}}{2\kappa (A + 2 c_{0}^{2})^{2}}\), and the constants \(C_f\), \(\tau\), \(A\), and \(B\) are defined as before.
\end{thm}
\begin{proof}
  The estimate of the first term \(E_{1}(x)\) has been given in \cref{eq:E1}. It remains to estimate the second term \(E_{2}(x)\) for the specific choices of \(n\), \(p_x\) and \(\gamma_x\). By the choices of \(n\) in \cref{eq:admissible_n}, \(p_x\) and \(\gamma_x\) in \cref{eq:gammax_px}, we have \(\frac{n^{2}\gamma_{x}^{2}}{2m^{2}} = \theta p_{x}\) with \(\theta>1\). Substituting it into the estimate of \(E_{2}(x)\) in \Cref{prop:E2} and applying Stirling's formula, we obtain
  \begin{align*}
    E_{2}(x) \leq \frac{A+2B}{A} (2n+1) \frac{\se^{p_{x}}}{\sqrt{2\pi p_{x}} p_{x}^{p_{x}}} \se^{-\theta p_{x}} (\theta p_{x})^{p_{x}} = \frac{A+2B}{A} \frac{2n+1}{\sqrt{2\pi p_{x}}} \se^{-(\theta - 1 - \ln(\theta))p_{x}}.
  \end{align*}
  The desired result then follows from substituting the choice of \(n\) in \cref{eq:admissible_n} and \(p_x\) in \cref{eq:gammax_px} into the above inequality.

\end{proof}

We remark that the conditions on \(\alpha\) and \(\kappa\) in \Cref{thm:filter_reconstruction} ensure that \(\theta>1\), which implies \(\theta-1-\ln(\theta)>0\). The second term in the error estimate in \Cref{thm:filter_reconstruction} then decays exponentially as \(m\) increases.

However, the filter reconstruction \(f_{m}^{\filter}\) could only achieve high accuracy when \(x\) is away from the jump discontinuity \(\xi\). Indeed, when \(d(x)\) is small, the error estimate in \Cref{thm:filter_reconstruction} could be very large. In particular, when \(x\) is close to \(\xi\), the error estimate could be even larger than one. Thus, we will employ a stable extrapolation method to recover the function values near the jump discontinuity in the next section.

\section{Hybrid filter-extrapolation method} \label{sec:hybrid}
In this section, we will combine the filter method and a stable extrapolation method to reconstruct the underlying piece-wise smooth function \(f\) in the whole region \([0,1]\) from its non-uniform Fourier data. Specifically, we will use the filter reconstruction \(f_{m}^{\filter}\) to recover the function values away from the jump discontinuity \(\xi\) and employ a stable extrapolation method to recover the function values near \(\xi\). We will show that the proposed hybrid method could achieve uniform exponential accuracy through the whole region \([0,1]\).

We will employ the stable extrapolation method proposed in \cite{Demanet2018} to extrapolate the function values near the jump discontinuity \(\xi\) based on the function values away from \(\xi\) recovered by the filter reconstruction \(f_{m}^{\filter}\). To this end, we first review the stable extrapolation method in \cite{Demanet2018}. Specifically, suppose a function \(g\) is analytically continuable to a function that is analytic in the open ellipse with foci \(\pm 1\), semimajor axis and semiminor axis lengths summing to \(\rho\), and bounded by \(Q>0\) there. We would like to extrapolate the function \(g\) on the interval \([1,\frac{\rho + \rho^{-1}}{2})\) based on its noisy samples \(\tilde{g}(x_{j}) = g(x_{j}) + \epsilon_{j}\) at the points \(x_{j} = \frac{2j}{N}-1\) for \(0\leq j\leq N\), where \(\epsilon_{j}\) is the noise with \(|\epsilon_{j}|\leq \epsilon\). The stable extrapolation method in \cite{Demanet2018} first computes the least squares polynomial approximation of degree \(M\):
\begin{align*}
  q_{M}=\arg\min_{q\in \cP_{M}} \sum_{j=0}^{N} |\tilde{g}(x_{j}) - q(x_{j})|^{2},
\end{align*}
where \(\cP_{M}\) is the space of all polynomials of degree at most \(M\). The solution to the above least squares problem is then used to extrapolate the function \(g\) on the interval \([1,\frac{\rho + \rho^{-1}}{2})\). It has been shown in \cite{Demanet2018} that if we choose \(M^{*}= \frac{\log \frac{Q}{\epsilon}}{\log \rho}\) and \(N\geq 4(M^{*})^{2}\) then the extrapolation could achieve the following exponential accuracy:
\begin{align}\label{eq:extrapolation_error}
  |g(x) - q_{M^{*}}(x)| \leq C_{\rho, \epsilon} \frac{Q^{1-\alpha(x)}}{1-r(x)} \epsilon^{\alpha(x)}, \quad x\in [1,\frac{\rho + \rho^{-1}}{2}),
\end{align}
where \(r(x) = \frac{x + \sqrt{x^{2}-1}}{\rho}\in (0,1)\), \(\alpha(x) = -\frac{\log r(x)}{\log \rho}\in (0,1)\), and \(C_{\rho, \epsilon}\) is a constant depending on \(\rho\) and polylogarithmically on \(\frac{1}{\epsilon}\). We could use Chebyshev polynomials to solve the above least squares problem efficiently and stably; see \cite{Demanet2018} for more details.

For the piece-wise smooth function \(f\) supported on \([0,1]\) with a single jump discontinuity \(\xi\in (0,1)\), we will use the filter reconstruction \(f_{m}^{\filter}\) to recover the function values on the intervals \([0,\xi - \delta]\) and \([\xi + \delta, 1]\) for some small \(\delta>0\) and then employ the stable extrapolation method to recover the function values on the interval \((\xi - \delta, \xi + \delta)\). In particular, we will extrapolate the function values on the interval \((\xi - \delta, \xi)\) based on the function values of \(f_{m}^{\filter}\) on the interval \([0, \xi - \delta]\) and extrapolate the function values on the interval \((\xi, \xi + \delta)\) based on the function values of \(f_{m}^{\filter}\) on the interval \([\xi + \delta, 1]\). We will show that the proposed hybrid method could achieve uniform exponential accuracy through the whole region \([0,1]\).

In the simplicity of presentation, we will only present the details of the stable extrapolation on the interval \((\xi - \delta, \xi)\) based on the function values of \(f_{m}^{\filter}\) on the interval \([0, \xi - \delta]\). The stable extrapolation on the interval \((\xi, \xi + \delta)\) could be done in a similar way. Specifically, for \(\delta > 0\) and \(M\in \bN\), we will define the extrapolated function \(f_{M, \delta}^{\extrap}\) on the interval \((\xi - \delta, \xi)\) as follows:
\begin{align}\label{eq:extrapolated_function}
  f_{M, \delta}^{\extrap}(x) = \arg\mathop{\min}_{q\in \cP_{M}} \sum_{j=0}^{N} |f_{m}^{\filter}(x_{j}) - q(x_{j})|^{2}, \quad x\in (\xi - \delta, \xi),
\end{align}
where \(x_{j} = \frac{(\xi - \delta)j}{N}\) for \(0\leq j\leq N\). When choosing the parameters \(\delta\), \(M\), and \(N\) appropriately, we have the following convergence result of the extrapolated function \(f_{M, \delta}^{\extrap}\).

\begin{prop}\label{prop:extrapolation}
  Suppose the underlying function \(f\) restricted on the interval \([0,\xi)\) is analytically continuable to a function that is analytic in the open ellipse with foci \(0\) and \(\xi\), semimajor axis and semiminor axis lengths summing to \(\rho\) for some \(\rho>\frac{\xi}{2}\) and bounded by \(Q>0\) there. For any \(0<\delta<\xi\) and \(\epsilon_{\delta}\geq \mathop{\max}_{x\in [0, \xi - \delta]} \left\vert f(x) - f_{m}^{\filter}(x) \right\vert \), if we choose \(M = \frac{\log \frac{Q}{\epsilon_{\delta}}}{\log \rho}\) and \(N\geq 4M^{2}\), then we have
  \begin{align*}
    |f(x) - f_{M, \delta}^{\extrap}(x)| \leq C_{\rho, \epsilon_{\delta}} \frac{Q^{1-\alpha_{\delta}(x)}}{1-r_{\delta}(x)} \epsilon_{\delta}^{\alpha_{\delta}(x)}, \quad x\in (\xi - \delta, \xi),
  \end{align*}
  where \(r_{\delta}(x) = \frac{(2x - \xi + \delta) + \sqrt{(2x - \xi + \delta)^{2} - (\xi-\delta)^{2}}}{2\rho}\in (0,1)\), \(\alpha_{\delta}(x) = -\frac{\log r_{\delta}(x)}{\log \rho - \log \frac{\xi-\delta}{2}}\in (0,1)\), and \(C_{\rho, \epsilon_{\delta}}\) is a constant depending on \(\rho\) and polylogarithmically on \(\frac{1}{\epsilon_{\delta}}\).
\end{prop}
\begin{proof}
  We first make a change of variable \(t = \frac{2x}{\xi - \delta} - 1\) to map the interval \([0, \xi - \delta]\) to the interval \([-1,1]\). Then the function \(g(t) = f\left(\frac{\xi - \delta}{2} (t+1)\right)\) is analytically continuable to a function that is analytic in the open ellipse with foci \(\pm 1\), semimajor axis and semiminor axis lengths summing to \(\frac{2\rho}{\xi - \delta}\) and bounded by \(Q>0\) there. The noisy samples of the function \(g\) at the points \(t_{j} = \frac{2j}{N}-1\) for \(0\leq j\leq N\) are given by
  \begin{align*}
    \tilde{g}(t_{j}) = f_{m}^{\filter}(x_{j}),
  \end{align*}
  where \(x_{j} = \frac{(\xi - \delta)j}{N}\) for \(0\leq j\leq N\), and the noise satisfies
  \begin{align*}
    |g(t_{j}) - \tilde{g}(t_{j})| = |f(x_{j}) - f_{m}^{\filter}(x_{j})| \leq \epsilon_{\delta}.
  \end{align*}
  Thus, we can view the samples of the function \(f_{m}^{\filter}\) on the interval \([0, \xi - \delta]\) as noisy samples of the function \(g\) with noise level at most \(\epsilon_{\delta}\). The extrapolated function \(f_{M, \delta}^{\extrap}\) on the interval \((\xi - \delta, \xi)\) could then be viewed as the least squares polynomial approximation \(q_{M}\) of degree \(M\) to the noisy samples \(\tilde{g}(t_{j})\) at the points \(t_{j} = \frac{2j}{N}-1\) for \(0\leq j\leq N\). The desired result then follows from applying the extrapolation error estimate \cref{eq:extrapolation_error} to the function \(g\) with noise level at most \(\epsilon_{\delta}\).
\end{proof}

We now present the hybrid filter-extrapolation method to reconstruct the underlying piece-wise smooth function \(f\) in the region \([0,\xi) \) from its non-uniform Fourier data. Specifically, for \(\delta > 0\), \(M\in \bN\), and \(N\in \bN\), we define the hybrid reconstruction \(f_{m, M, \delta}^{\hybrid}\) on the interval \([0,\xi)\) as follows:
    \begin{align*}
      f_{m, M, \delta}^{\hybrid}(x) = \begin{cases}
                                        f_{m}^{\filter}(x),         & x\in [0, \xi - \delta],   \\
                                        f_{M, \delta}^{\extrap}(x), & x\in (\xi - \delta, \xi).
                                      \end{cases}
    \end{align*}

    We have the following convergence result of the hybrid reconstruction \(f_{m, M, \delta}^{\hybrid}\).
    \begin{thm}\label{thm:hybrid_reconstruction}
      Suppose the conditions in \Cref{thm:filter_reconstruction} and \Cref{prop:extrapolation} hold. For any \(0<\delta<\xi\), there exists a positive constant \(C_{\delta}\) such that
      \begin{align*}
        \left\vert f(x) - f_{m, M, \delta}^{\hybrid}(x) \right\vert \leq C_{\delta} m^{\frac{9}{4}}\se^{-\eta m \delta}, \quad x\in [0, \xi - \delta],
      \end{align*}
      where \(\eta = \min\left\{\ln \tau, (\theta - 1 - \ln(\theta))\kappa\right\}\) with \(\tau\) and \(\theta\) defined in \Cref{thm:filter_reconstruction}. Let \(\epsilon^{\delta} = C_{\delta} m^{\frac{9}{4}}\se^{-\eta m \delta}, \quad x\in [0, \xi - \delta]\) and choose \(M = \frac{\log \frac{Q}{\epsilon_{\delta}}}{\log \rho}\) and \(N = 4M^{2}\). Then we have
      \begin{align*}
        |f(x) - f_{m, M, \delta}^{\hybrid}(x)| \leq C_{\rho, \epsilon_{\delta}} \frac{Q^{1-\alpha_{\delta}^{*}}}{1-r_{\delta}^{*}} (\epsilon_{\delta})^{\alpha_{\delta}^{*}}, \quad x\in (\xi - \delta, \xi),
      \end{align*}
      where \(r_{\delta}^{*} = \frac{\xi + \delta + 2\sqrt{\xi\delta}}{2\rho}\) and \(\alpha_{\delta}^{*} = -\frac{\log r_{\delta}^{*}}{\log \rho - \log \frac{\xi-\delta}{2}}\).
    \end{thm}
    \begin{proof}
      The estimate on the interval \([0, \xi - \delta]\) follows directly from the error estimate in \Cref{thm:filter_reconstruction} with \(d(x)\geq \delta\) for any \(x\in [0, \xi - \delta]\). For the estimate on the interval \((\xi - \delta, \xi)\), it follows from applying the error estimate in \Cref{prop:extrapolation} with \(\epsilon_{\delta} = C_{\delta} m^{\frac{9}{4}}\se^{-\eta m \delta}\) and noting that \(r_{\delta}(x)\) is an increasing function of \(x\) on the interval \((\xi - \delta, \xi)\).
    \end{proof}

    We remark that the error estimate in \Cref{thm:hybrid_reconstruction} shows that the proposed hybrid filter-extrapolation method could achieve uniform exponential accuracy through the whole region \([0,\xi)\). It does not depend on how close \(x\) is to the jump discontinuity \(\xi\). The same result could be obtained for the hybrid reconstruction on the interval \((\xi, 1]\) similarly. Thus, the proposed hybrid method could achieve uniform exponential accuracy through the whole region \([0,1]\). Additionally, we note that the parameter \(\alpha_{\delta}^{*}\) in the error estimate on the interval \((\xi - \delta, \xi)\) satisfies \(0<\alpha_{\delta}^{*}<1\) and the extrapolation error on the interval \((\xi - \delta, \xi)\) is always worse than the approximation error \(\epsilon_{\delta}\) on the interval \([0, \xi - \delta]\). This suggests that we should choose \(\delta\) to minimize the extrapolation error to achieve the best overall accuracy. In other words, we will consider choosing the parameter \(\delta\) to minimize \((\epsilon_{\delta})^{\alpha_{\delta}^{*}}\), which is equivalent to
    \begin{align*}
      \mathop{\min}_{0<\delta<\xi} -\frac{\log \frac{\xi + \delta + 2\sqrt{\xi\delta}}{2\rho}}{\log \rho - \log \frac{\xi-\delta}{2}} \log \left(C_{\delta} m^{\frac{9}{4}}\se^{-\eta m \delta}\right).
    \end{align*}
    The optimal choice of \(\delta\) could be obtained by solving the above one-dimensional optimization problem numerically. Similar optimization could be done for the hybrid reconstruction on the interval \((\xi, 1]\).

When the underlying piece-wise smooth function \(f\) has multiple jump discontinuities, we can apply the proposed hybrid filter-extrapolation method in each sub-interval separated by the jump discontinuities to achieve uniform exponential accuracy through the whole region \([0,1]\).

\section{Numerical experiments} \label{sec:numerical}
In this section, we present some numerical experiments to demonstrate the performance of the proposed hybrid filter-extrapolation method for reconstructing a piece-wise smooth function \(f\) from its non-uniform Fourier data \(\hat{f}(\lambda_{j}), -m\leq j\leq m\).

Specifically, we will first compute the filter reconstruction \(f_{m}^{\filter}\) \cref{eq:filter_reconstruction}, which applies the admissible frame approximation to the mollified function \(\rho_{p_{x}, \gamma_{x}} * f\). The mollifier \(\rho_{p_{x}, \gamma_{x}}\) is defined in \cref{eq:mollifier} with the parameters \(p_{x}\) and \(\gamma_{x}\) chosen as in \cref{eq:gammax_px} with \(\alpha = 1\) and \(\kappa = \frac{1}{15}\) as suggested in \cite{Tanner2006}.
The admissible frame is the Fourier basis \(\{\phi_{j}(x) = \se^{2\pi \si j x}\}_{-n\leq j \leq n}\).

Suppose the jump discontinuities of the underlying piece-wise smooth function \(f\) divide the interval \([0,1]\) into subintervals.
On each subinterval \([\xi_{l}, \xi_{l+1}]\), we will employ the stable extrapolation method to recover the function values near the jump discontinuities \(\xi_{l}\) and \(\xi_{l+1}\) based on the function values of \(f_{m}^{\filter}\) away from the jump discontinuities.
Specifically, for \(\delta > 0\), \(M\in \bN\), and \(N\in \bN\), we will find the extrapolated function \(f_{M, \delta}^{\extrap}\) based on \(M\) equally spaced points on  \([\xi_{l} + \delta, \xi_{l+1} - \delta]\) as shown in \cref{eq:extrapolated_function}. We then use it to recover the function values on the intervals \((\xi_{l}, \xi_{l} + \delta)\) and \((\xi_{l+1} - \delta, \xi_{l+1})\).

The hybrid reconstruction \(f_{m, M, \delta}^{\hybrid}\) on the interval \([\xi_{l}, \xi_{l+1}]\) is then obtained by combining the filter reconstruction \(f_{m}^{\filter}\) on \([\xi_{l}+ \delta, \xi_{l+1} - \delta]\) and the extrapolated function \(f_{M, \delta}^{\extrap}\) on \((\xi_{l}, \xi_{l} + \delta)\) and \((\xi_{l+1} - \delta, \xi_{l+1})\). The overall hybrid reconstruction \(f_{m, M, \delta}^{\hybrid}\) on the interval \([0,1]\) is then obtained by combining the hybrid reconstructions on all subintervals.

We will consider two types of non-uniform Fourier frequencies \(\lambda_{j}\): jittered sampling and log sampling. The jittered sampling is defined as
\begin{align*}
  \lambda_{j} = j + \epsilon_{j}, \quad -m \leq j \leq m,
\end{align*}
where \(\epsilon_{j}\) is a random variable uniformly distributed in the interval \([-\frac{1}{4}, \frac{1}{4}]\). The number of admissible Fourier modes \(n\)  for the jittered sampling is chosen as \(n = 0.6 m\).

The log sampling is defined as
\begin{align*}
  \lambda_{j} =
  \begin{cases}
    -\exp (-v + \frac{v+\log m}{m-1} (|j| - 1)), & \text{ if } -m\leq j \leq  -1; \\
    0,                                           & \text{ if } j = 0;             \\
    \exp (-v + \frac{v+\log m}{m-1} (|j| - 1)),  & \text{ if }  1\leq j \leq m.
  \end{cases}
\end{align*}
where \(v=0.001\). The number of admissible Fourier modes \(n\) for the log sampling is chosen as \(n = 2m^{0.6}\).

\begin{example}
  We first consider the following piece-wise smooth function with a single jump discontinuity at \(\xi = 0.5\):
  \begin{align*}
    f_{1}(x) = \begin{cases}
                 \sin(4\pi x), & 0\leq x < 0.5,    \\
                 \sin(2\pi x), & 0.5\leq x \leq 1.
               \end{cases}
  \end{align*}

  We choose \(\delta = \frac{1}{40}\). For each \(m\in \{128, 256, 512\}\), we first compute the filter reconstruction \(f_{m}^{\filter}\). We then use its values on the intervals \([0, \xi - \delta]\) and \([\xi + \delta, 1]\) to set up the least squares problem \cref{eq:extrapolated_function} to find the extrapolated function \(f_{M, \delta}^{\extrap}\), where \(M = 4 + m\delta\) and \(N = 4M^{2}\). We combine \(f_{m}^{\filter}\) and \(f_{M, \delta}^{\extrap}\) to obtain the hybrid reconstruction \(f_{m, M, \delta}^{\hybrid}\).

  For jittered sampling, we display the filtered reconstruction \(f_{m}^{\filter}\) and the hybrid reconstruction \(f_{m, M, \delta}^{\hybrid}\) in \Cref{fig:jittered_reconstruction}. In addition, we display the pointwise approximation errors of the filter reconstruction \(f_{m}^{\filter}\) and the hybrid reconstruction \(f_{m, M, \delta}^{\hybrid}\) in \Cref{fig:jittered_error}.
  \begin{figure}[!h]
    \centering
    \begin{tabular}{ccc}
      \textbf{m=128}                                                                     & \textbf{m=256} & \textbf{m=512} \\
      \includegraphics[width=0.32\textwidth]{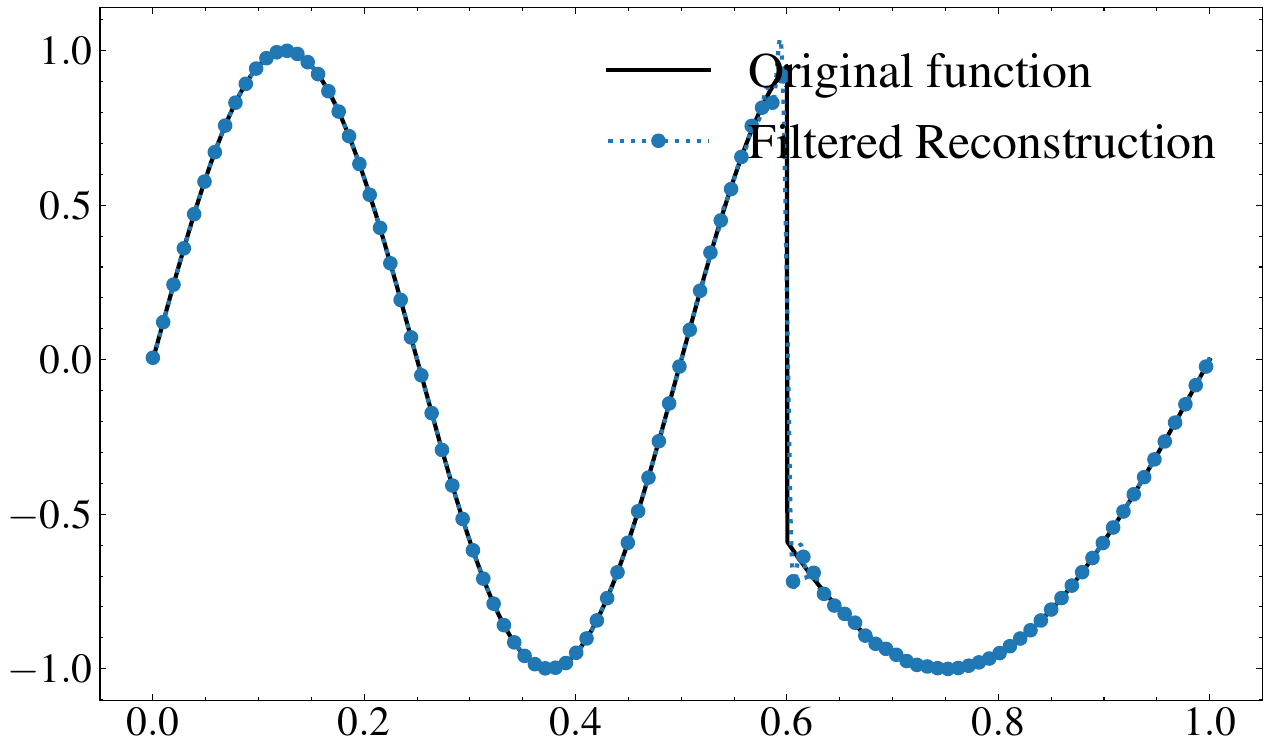} &
      \includegraphics[width=0.32\textwidth]{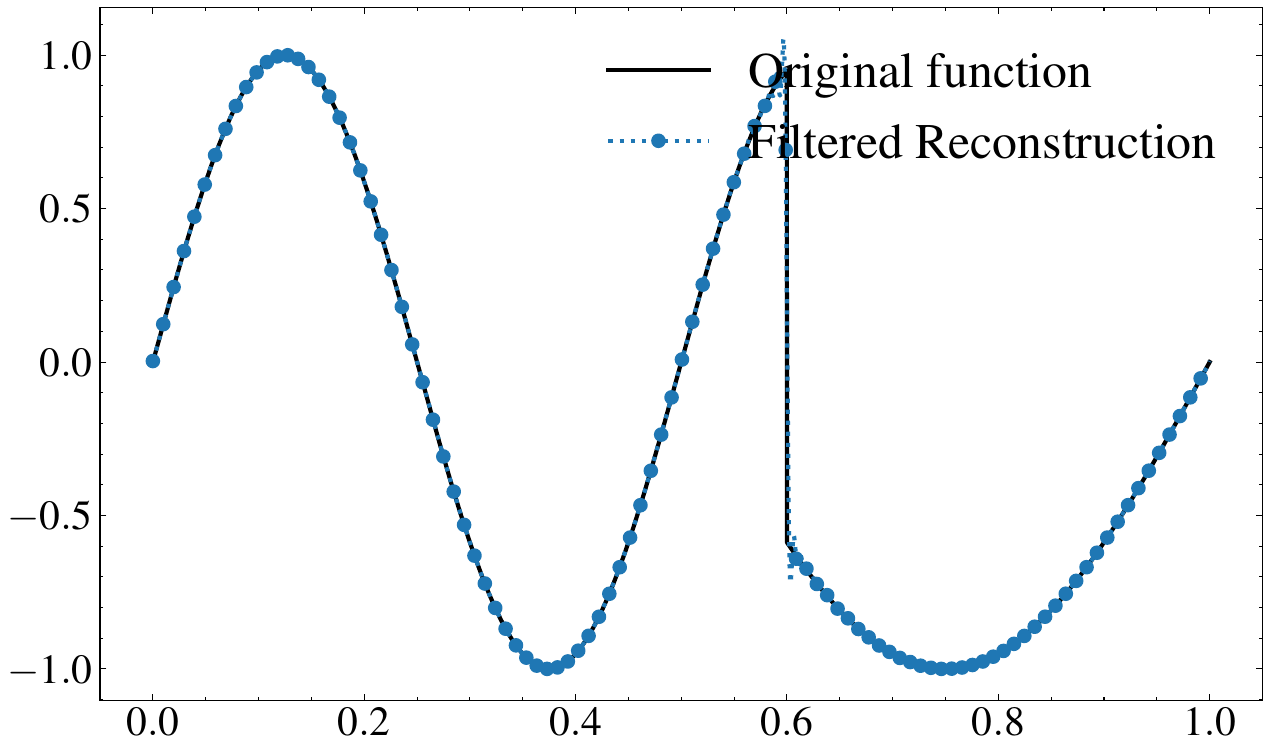} &
      \includegraphics[width=0.32\textwidth]{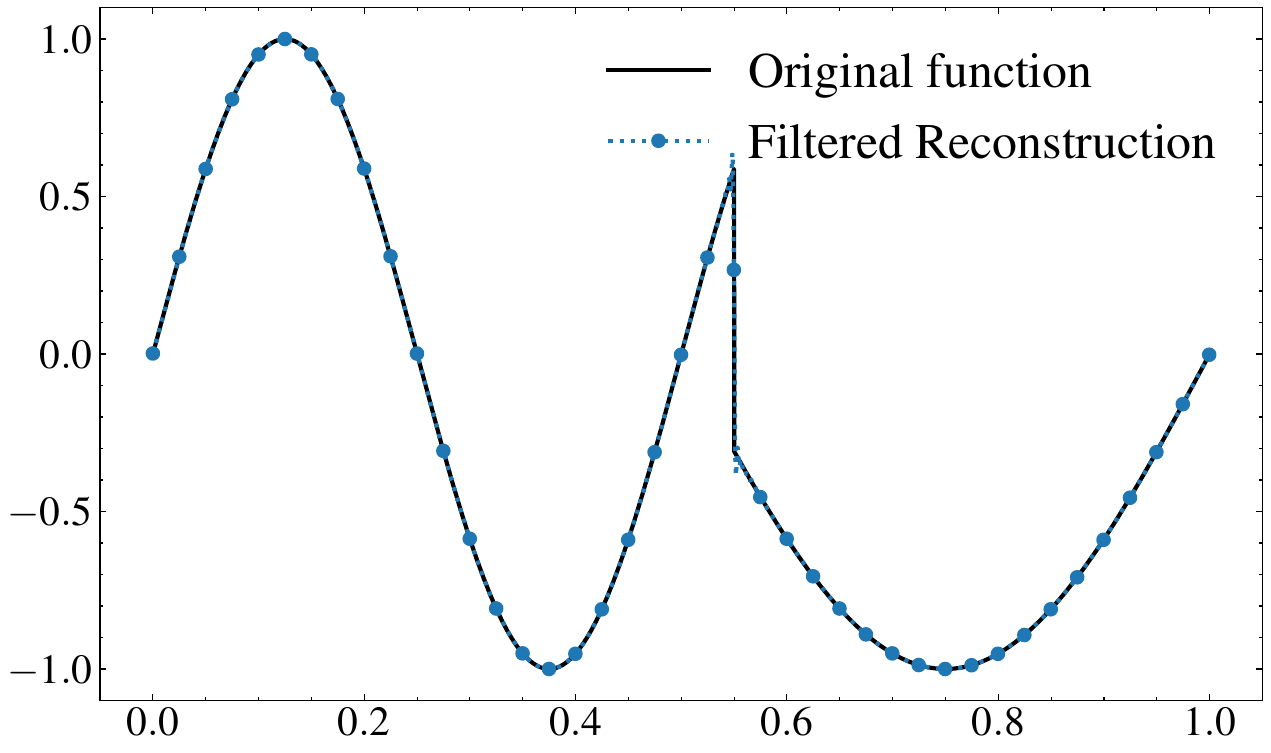}                                   \\
      \includegraphics[width=0.32\textwidth]{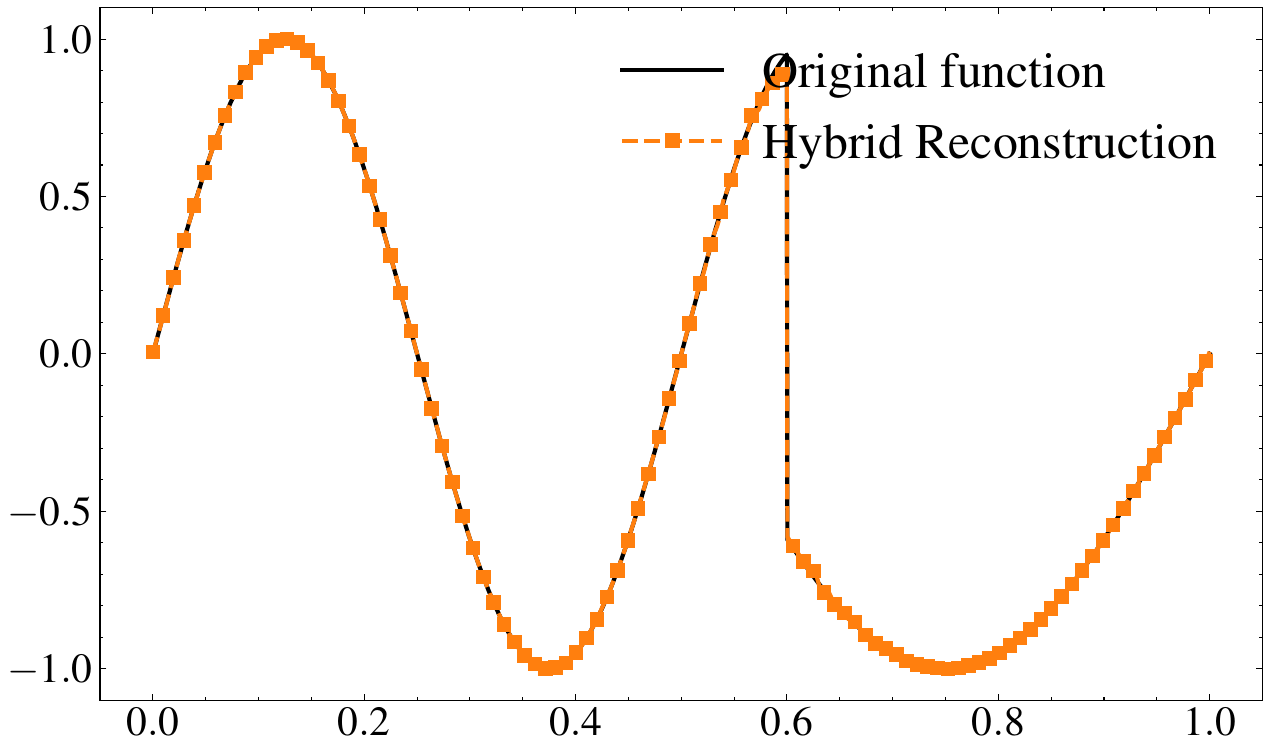}    &
      \includegraphics[width=0.32\textwidth]{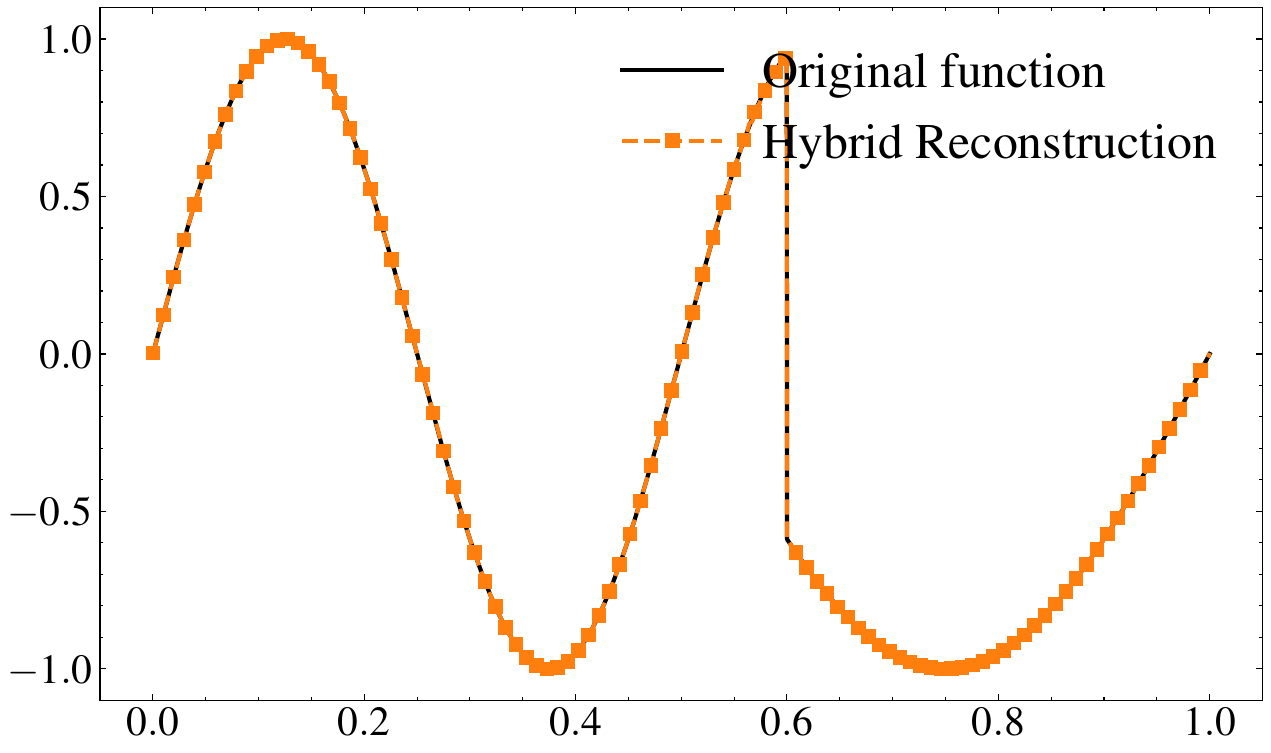}    &
      \includegraphics[width=0.32\textwidth]{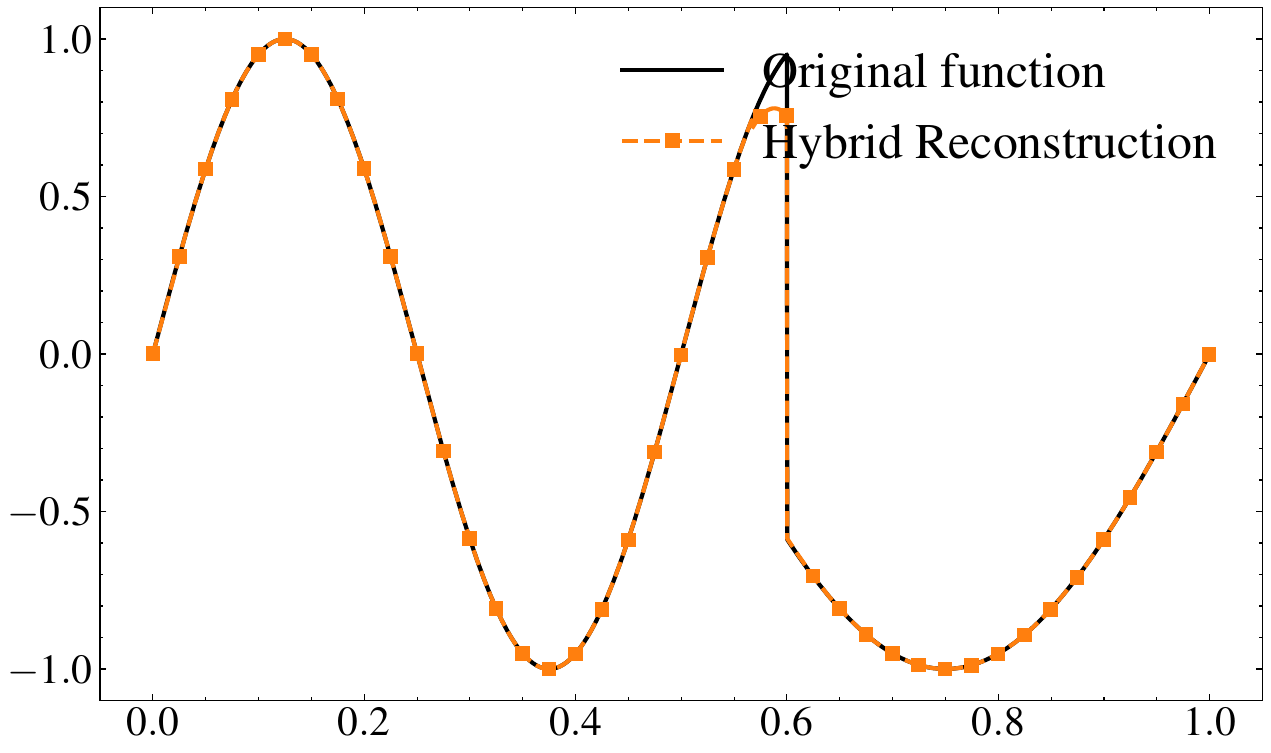}                                      \\
    \end{tabular}
    \caption{Filtered reconstruction \(f_{m}^{\filter}\) (top row) and hybrid reconstruction \(f_{m, M, \delta}^{\hybrid}\) (bottom row) for the jittered sampling with \(m=128, 256, 512\).}
    \label{fig:jittered_reconstruction}
  \end{figure}

  \begin{figure}[!h]
    \centering
    \begin{tabular}{ccc}
      \textbf{m=128}                                                                & \textbf{m=256} & \textbf{m=512} \\
      \includegraphics[width=0.32\textwidth]{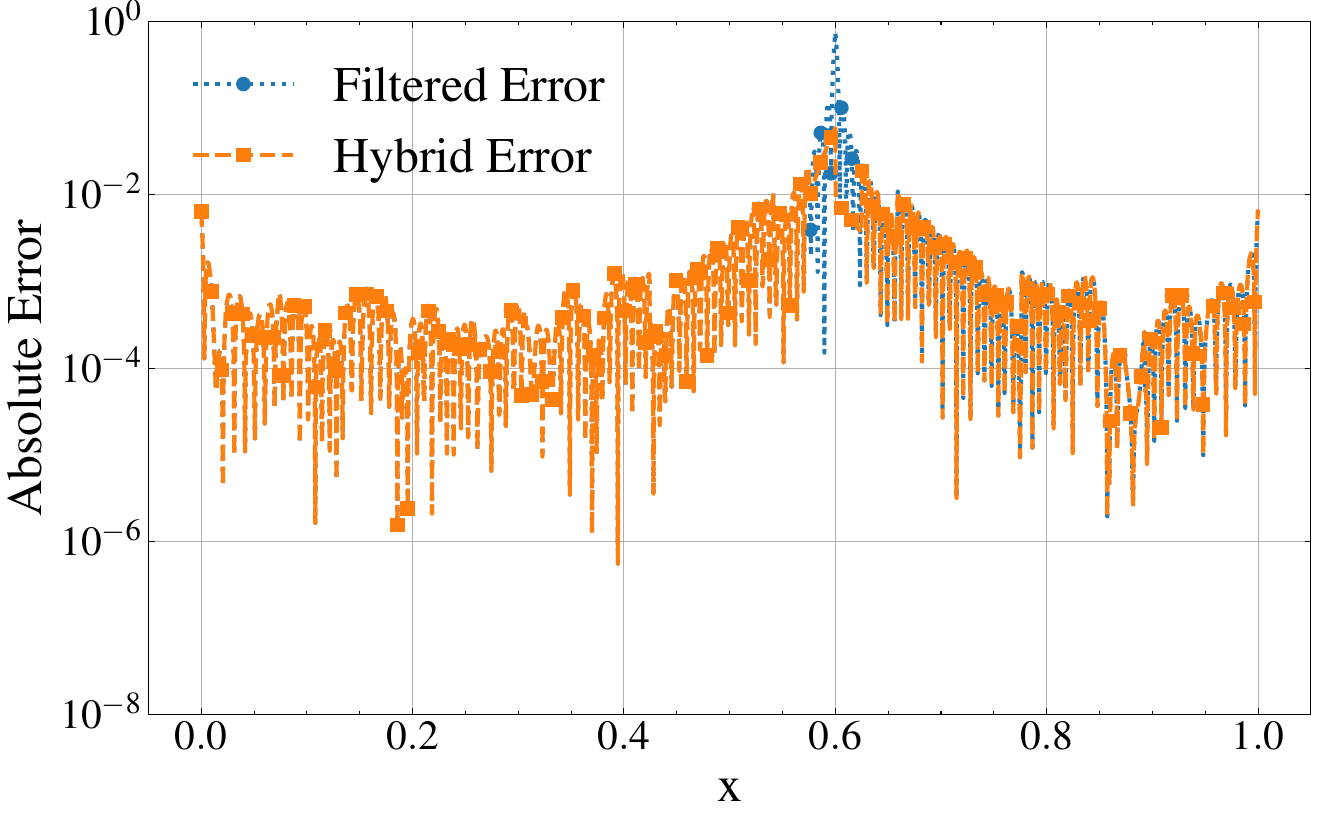} &
      \includegraphics[width=0.32\textwidth]{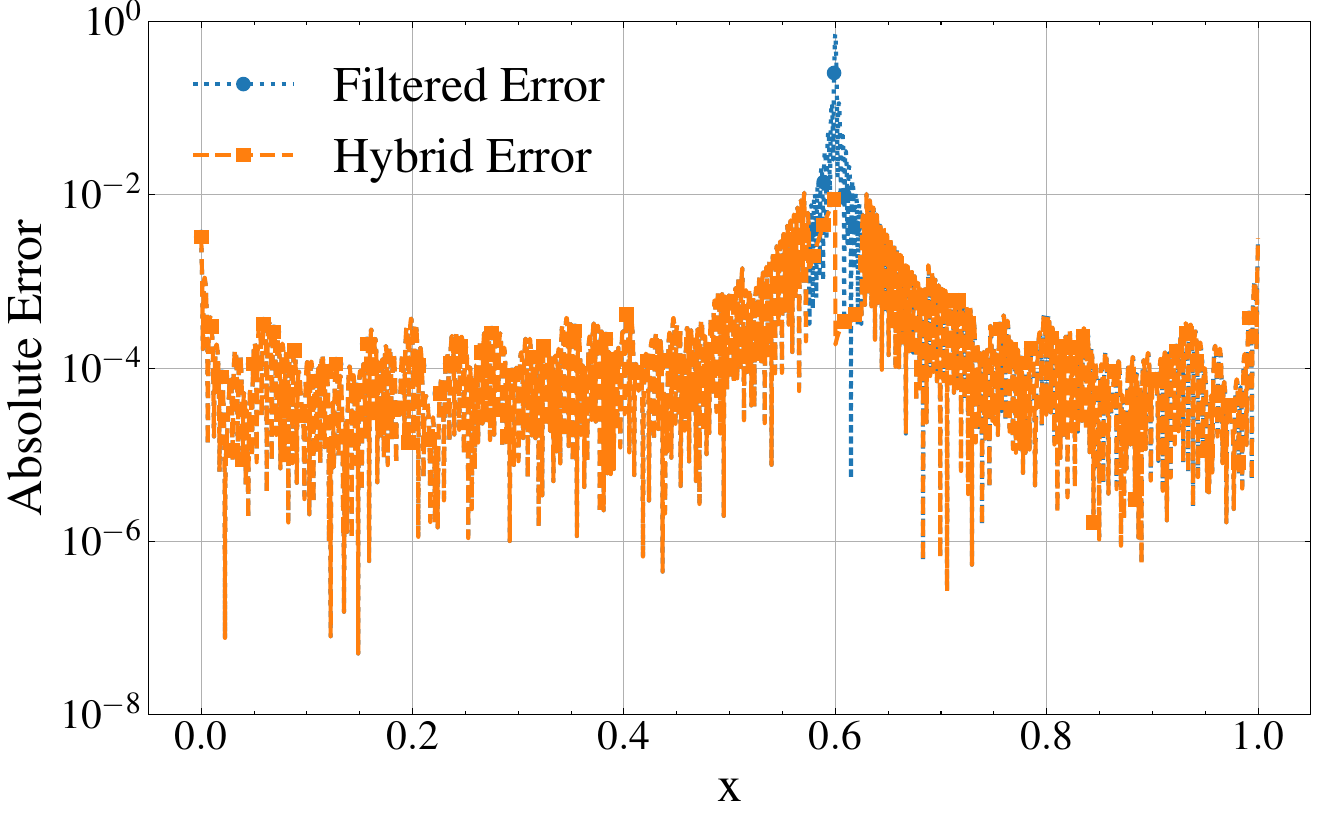} &
      \includegraphics[width=0.32\textwidth]{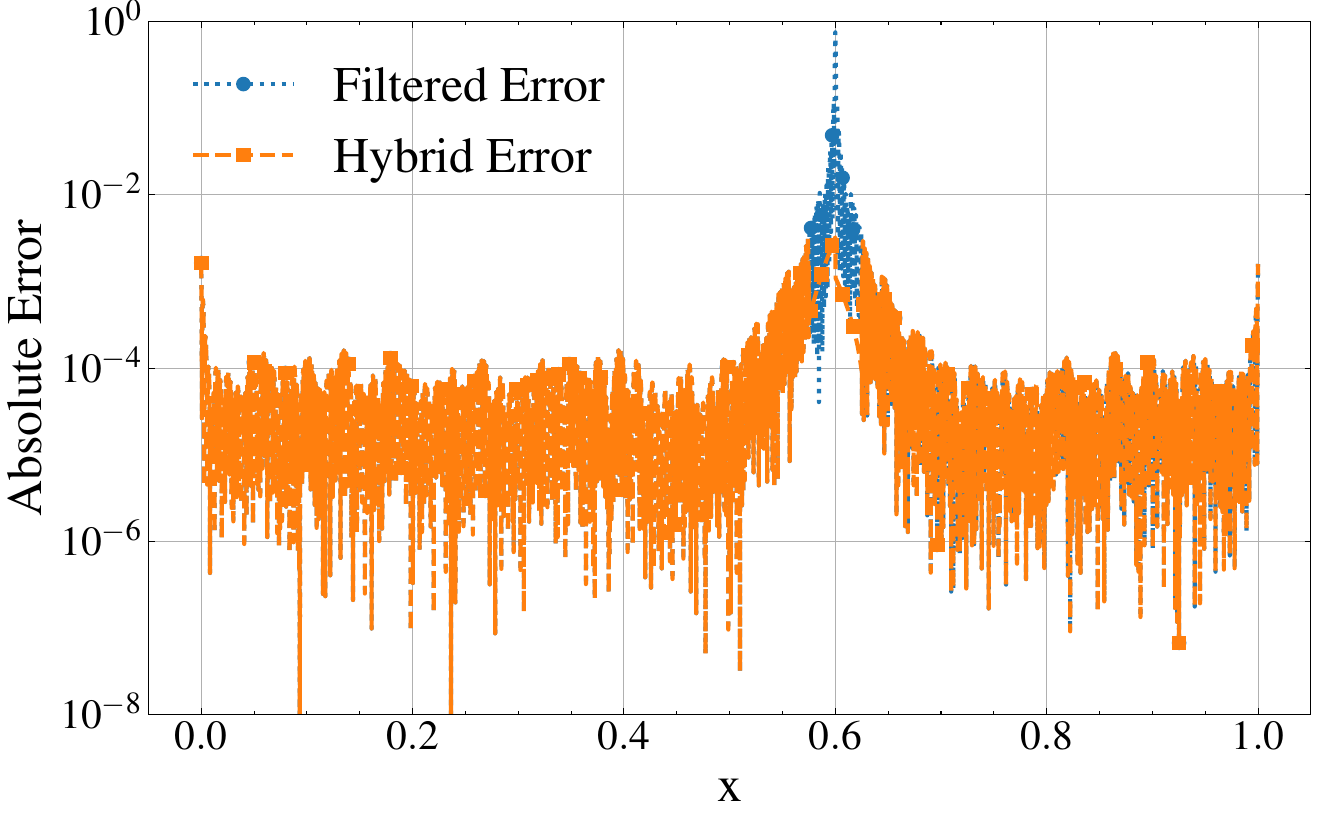}                                   \\
    \end{tabular}
    \caption{Pointwise approximation errors of the filter reconstruction \(f_{m}^{\filter}\)  and the hybrid reconstruction \(f_{m, M, \delta}^{\hybrid}\) for the jittered sampling with \(m=128, 256, 512\).}
    \label{fig:jittered_error}
  \end{figure}

  For log sampling, we display the reconstructed functions in \Cref{fig:log_reconstruction} and the pointwise approximation errors in \Cref{fig:log_error}.

  \begin{figure}[!h]
    \centering
    \begin{tabular}{ccc}
      \textbf{m=128}                                                                     & \textbf{m=256} & \textbf{m=512} \\
      \includegraphics[width=0.32\textwidth]{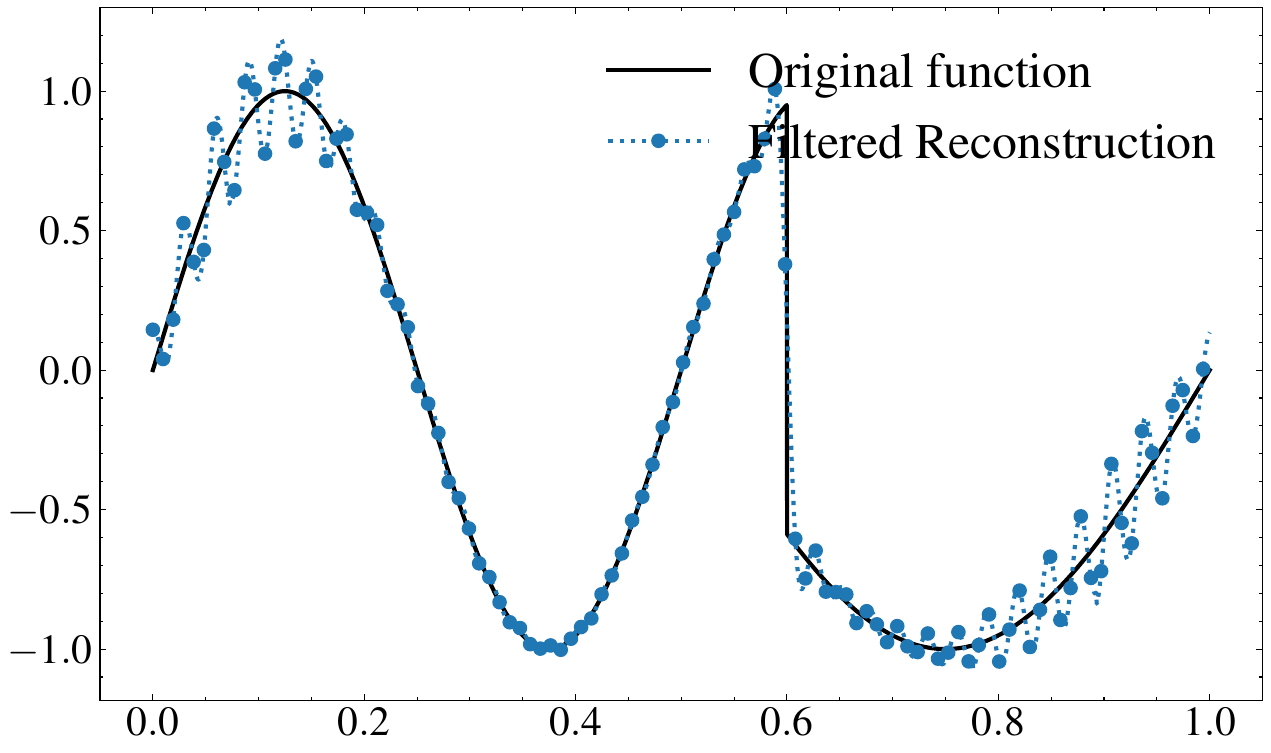} &
      \includegraphics[width=0.32\textwidth]{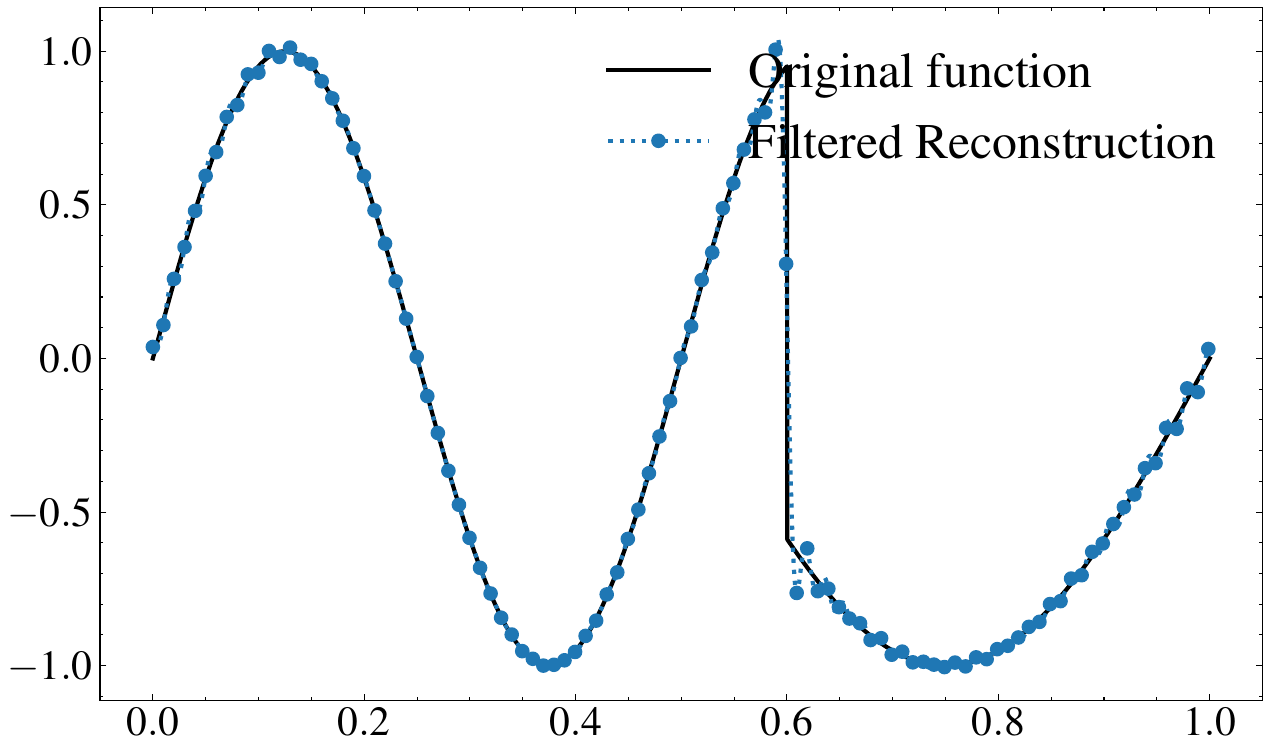} &
      \includegraphics[width=0.32\textwidth]{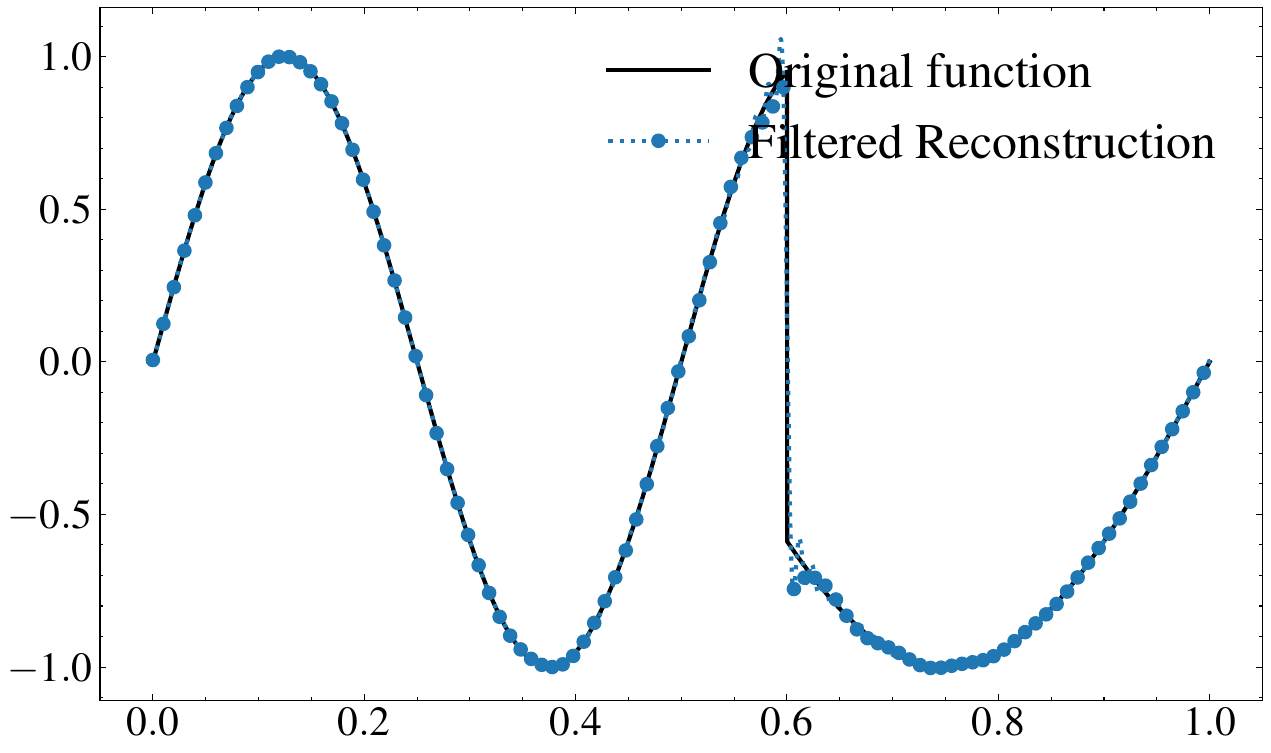}                                   \\
      \includegraphics[width=0.32\textwidth]{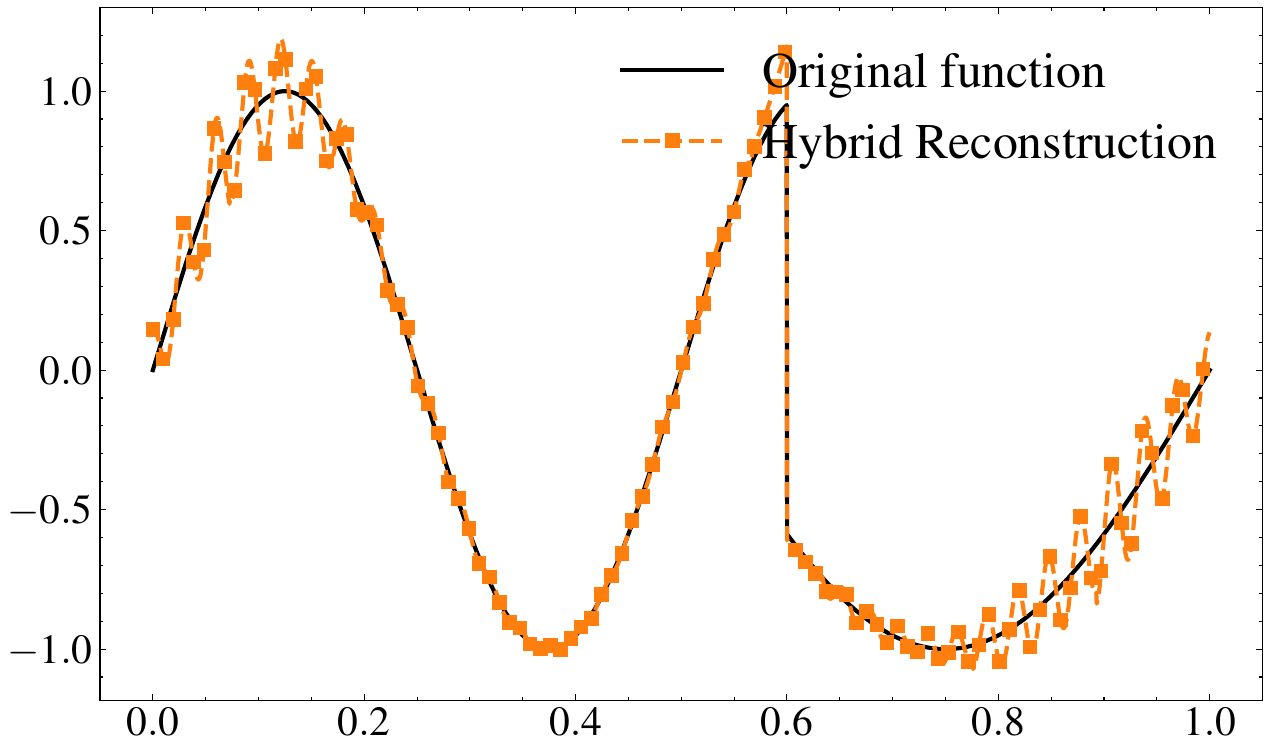}    &
      \includegraphics[width=0.32\textwidth]{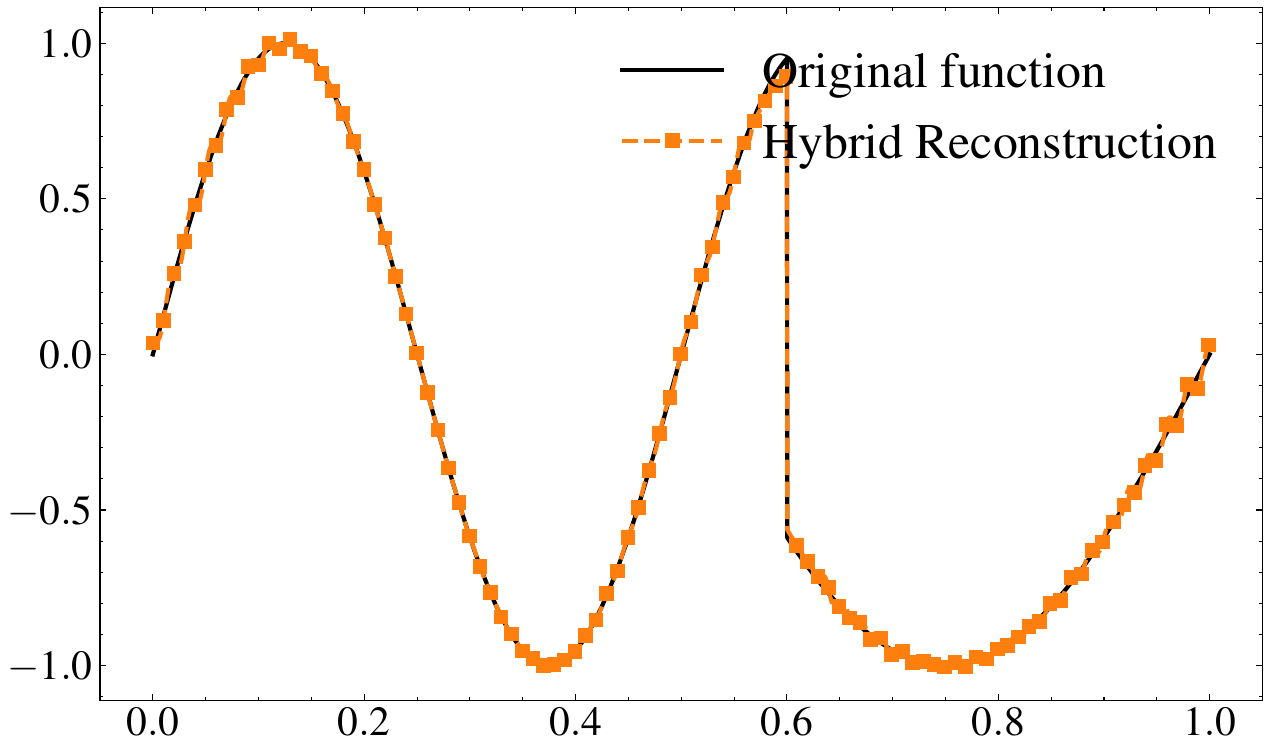}    &
      \includegraphics[width=0.32\textwidth]{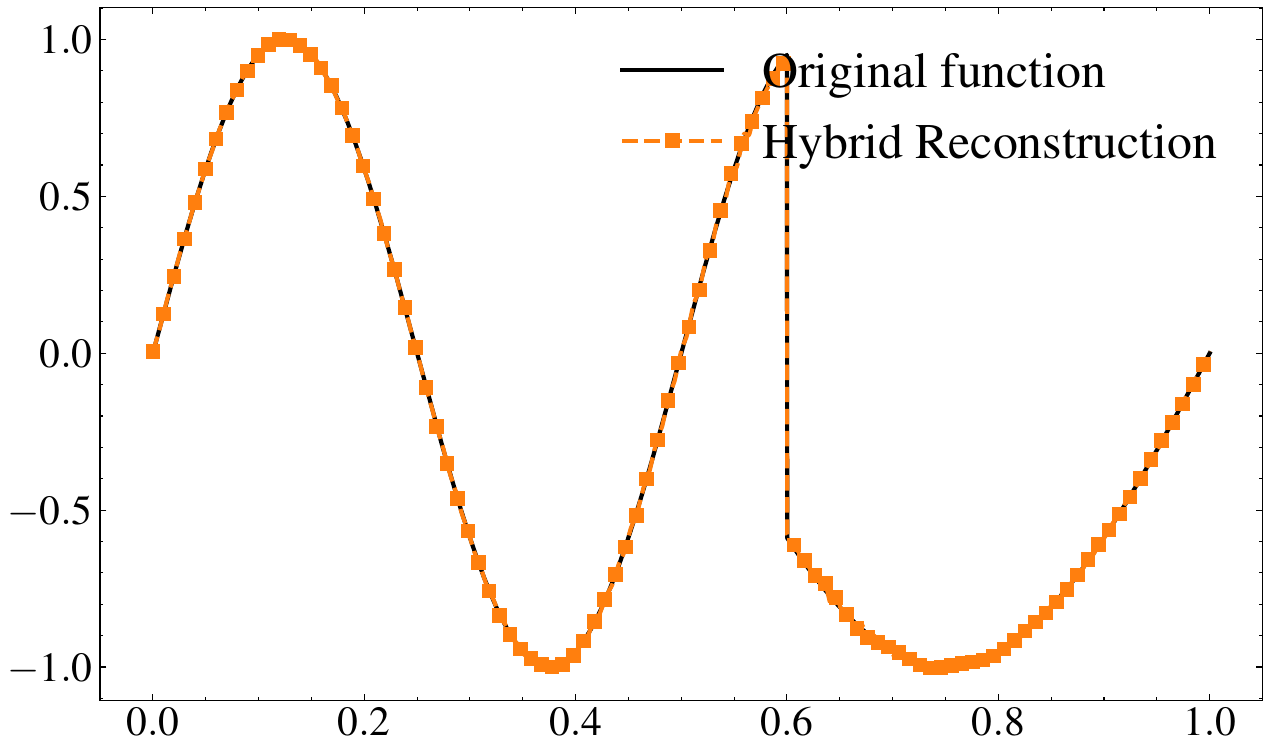}                                      \\
    \end{tabular}
    \caption{Filtered reconstruction \(f_{m}^{\filter}\) (top row) and hybrid reconstruction \(f_{m, M, \delta}^{\hybrid}\) (bottom row) for the log sampling with \(m=128, 256, 512\).}
    \label{fig:log_reconstruction}
  \end{figure}

  \begin{figure}[!h]
    \centering
    \begin{tabular}{ccc}
      \textbf{m=128}                                                                & \textbf{m=256} & \textbf{m=512} \\
      \includegraphics[width=0.32\textwidth]{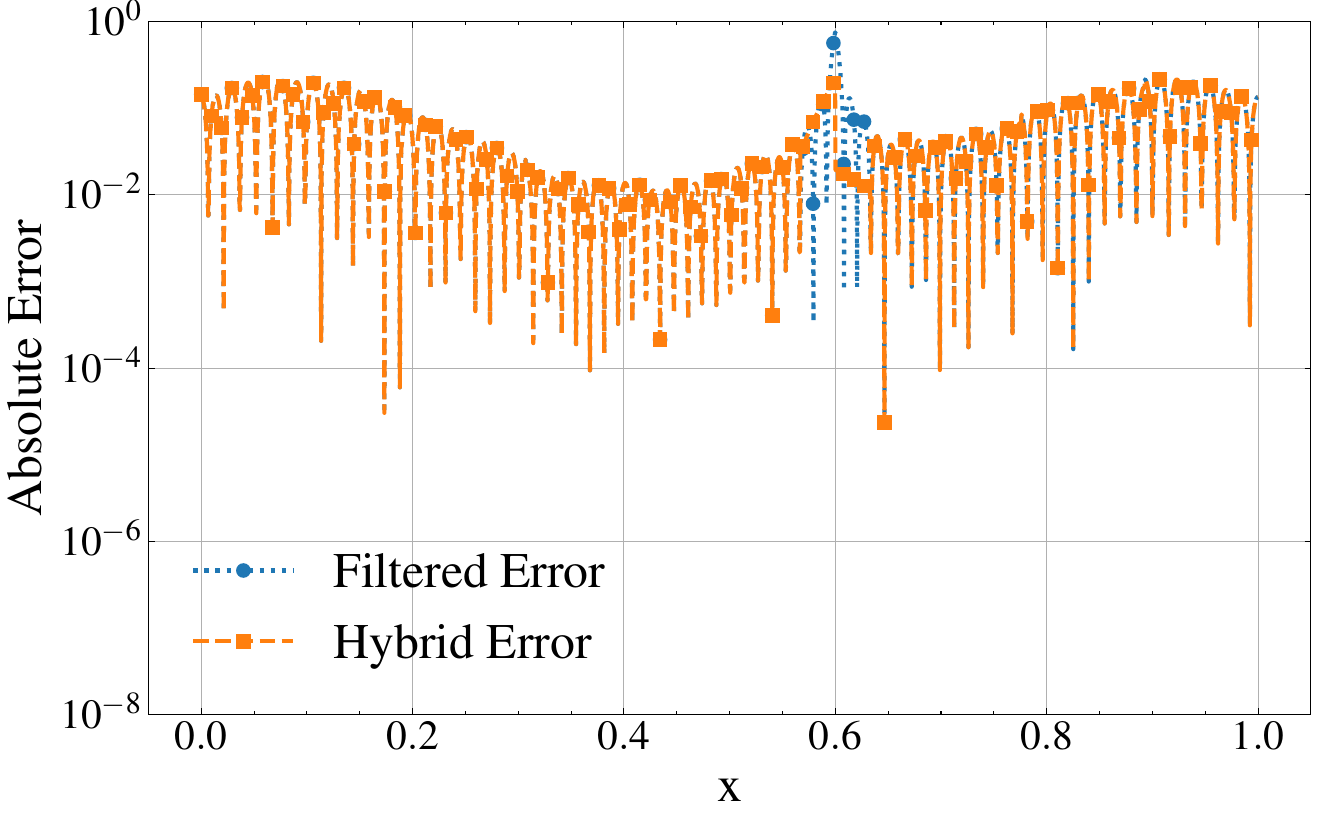} &
      \includegraphics[width=0.32\textwidth]{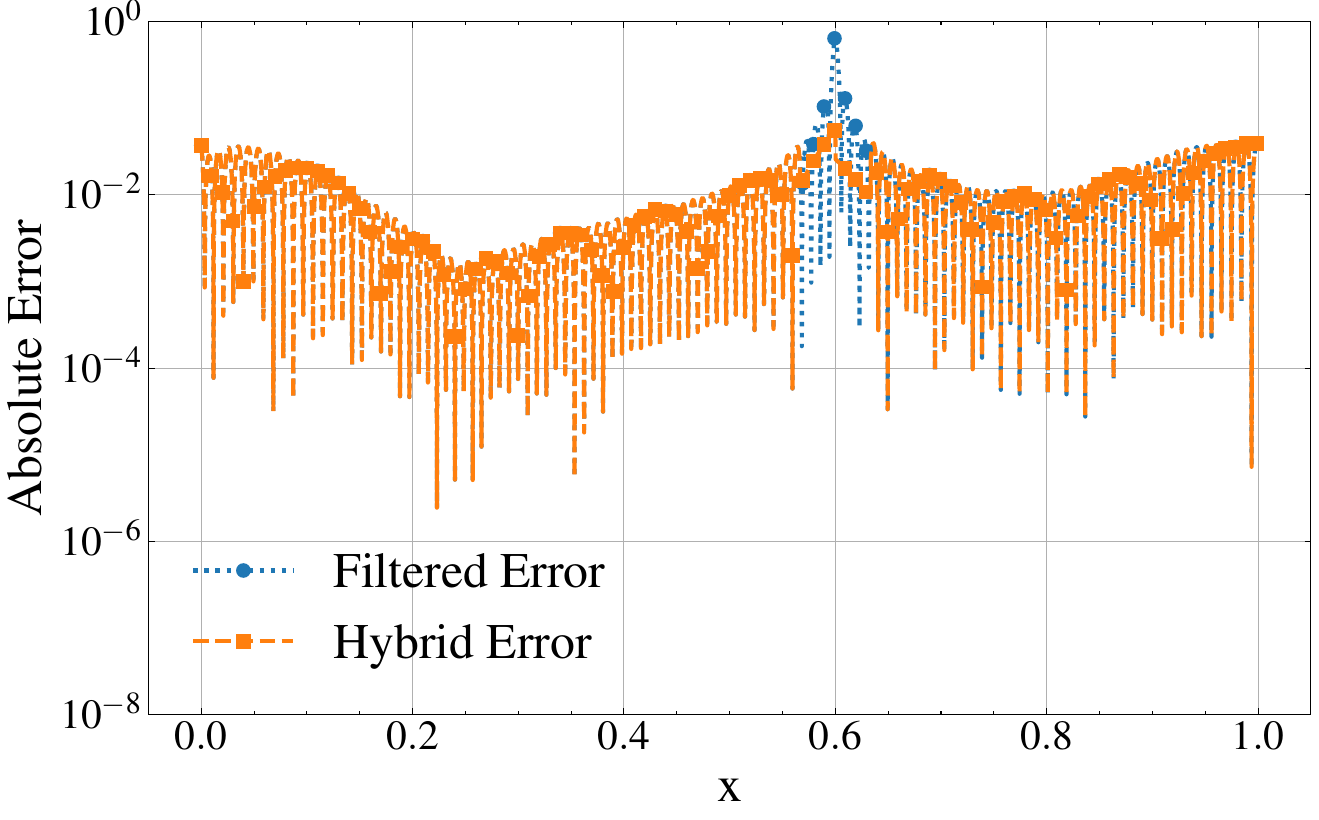} &
      \includegraphics[width=0.32\textwidth]{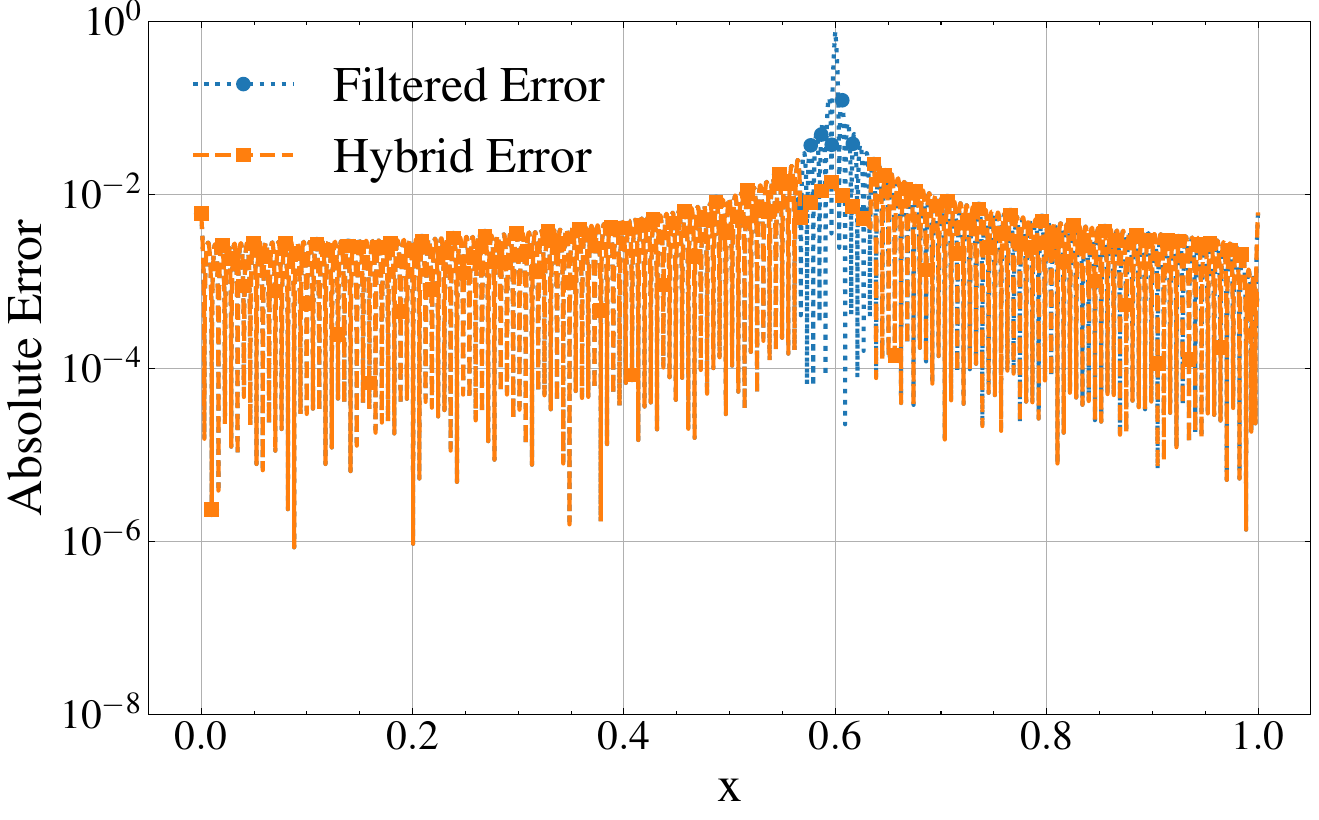}                                   \\
    \end{tabular}
    \caption{Pointwise approximation errors of the filter reconstruction \(f_{m}^{\filter}\)  and the hybrid reconstruction \(f_{m, M, \delta}^{\hybrid}\) for the log sampling with \(m=128, 256, 512\).}
    \label{fig:log_error}
  \end{figure}
\end{example}

\begin{example}
  We next consider the following piece-wise smooth function with multiple jump discontinuities at \(\xi_{1} = 0\), \(\xi_{2} = 0.3\), \(\xi_{3} = 0.6\), and \(\xi_{4} = 1\):
  \begin{align*}
    f_{2}(x) = \begin{cases}
                 \se^{-5x^{2}},       & 0\leq x < 0.3,   \\
                 \cos(2\pi x),        & 0.3\leq x < 0.7, \\
                 \se^{x}\sin(4\pi x), & 0.7\leq x < 1.
               \end{cases}
  \end{align*}

  We choose \(\delta = \frac{1}{40}\). Similarly, for each \(m\in \{128, 256, 512\}\), we first compute the filter reconstruction \(f_{m}^{\filter}\). We then use its values on the intervals \([\xi_{l} + \delta, \xi_{l+1} - \delta]\) for \(1\leq l\leq 3\) to set up the least squares problem \cref{eq:extrapolated_function} to find the extrapolated function \(f_{M, \delta}^{\extrap}\), where \(M = 4 + m\delta\) and \(N = 4M^{2}\). We combine \(f_{m}^{\filter}\) and \(f_{M, \delta}^{\extrap}\) to obtain the hybrid reconstruction \(f_{m, M, \delta}^{\hybrid}\).

  For jittered sampling, we display the filtered reconstruction \(f_{m}^{\filter}\) and the hybrid reconstruction \(f_{m, M, \delta}^{\hybrid}\) in \Cref{fig:jittered_reconstruction_multi} and the pointwise approximation errors of the filter reconstruction \(f_{m}^{\filter}\) and the hybrid reconstruction \(f_{m, M, \delta}^{\hybrid}\) in \Cref{fig:jittered_error_multi}.

  \begin{figure}[!h]
    \centering
    \begin{tabular}{ccc}
      \textbf{m=128}                                                                        & \textbf{m=256} & \textbf{m=512} \\
      \includegraphics[width=0.32\textwidth]{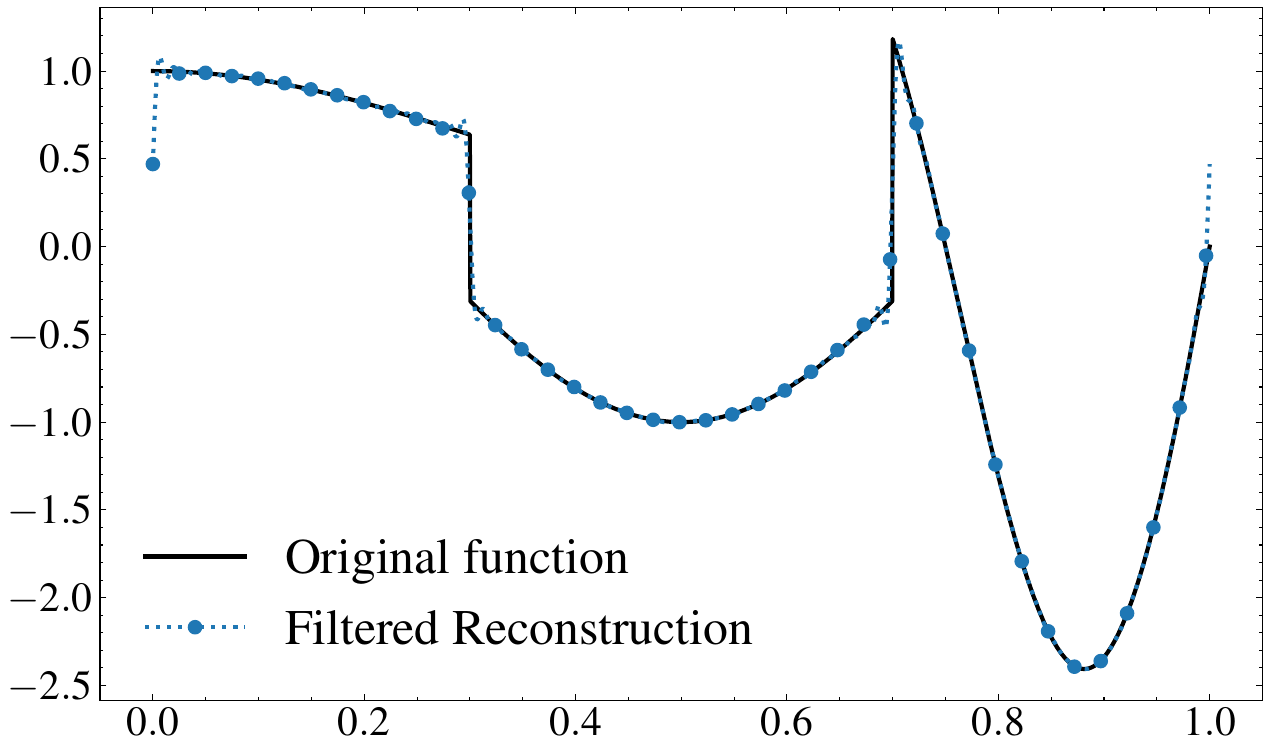} &
      \includegraphics[width=0.32\textwidth]{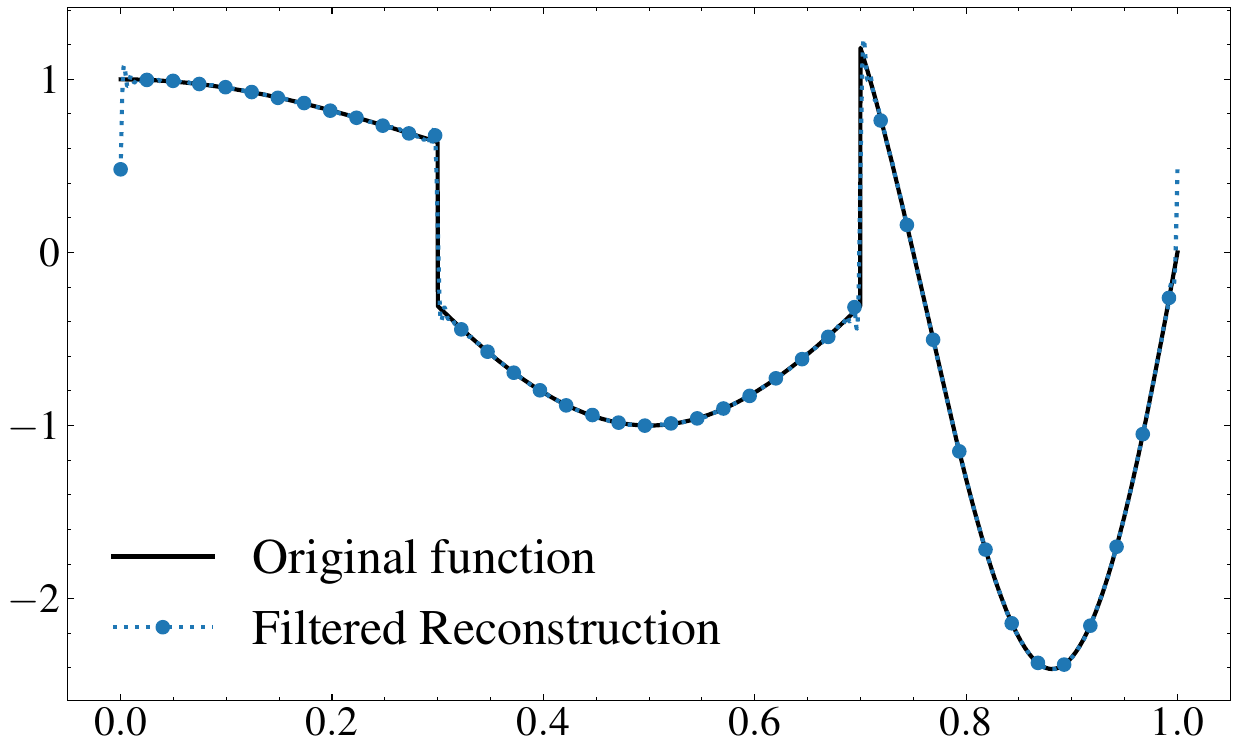} &
      \includegraphics[width=0.32\textwidth]{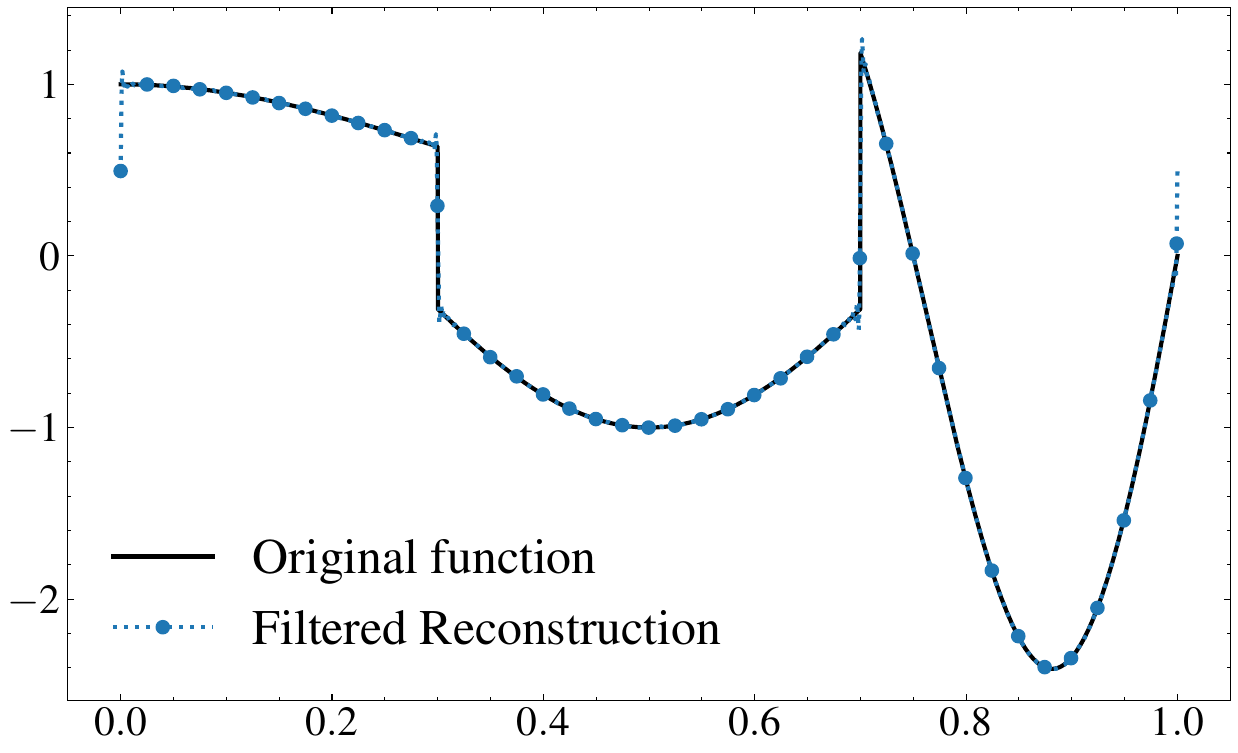}                                   \\
      \includegraphics[width=0.32\textwidth]{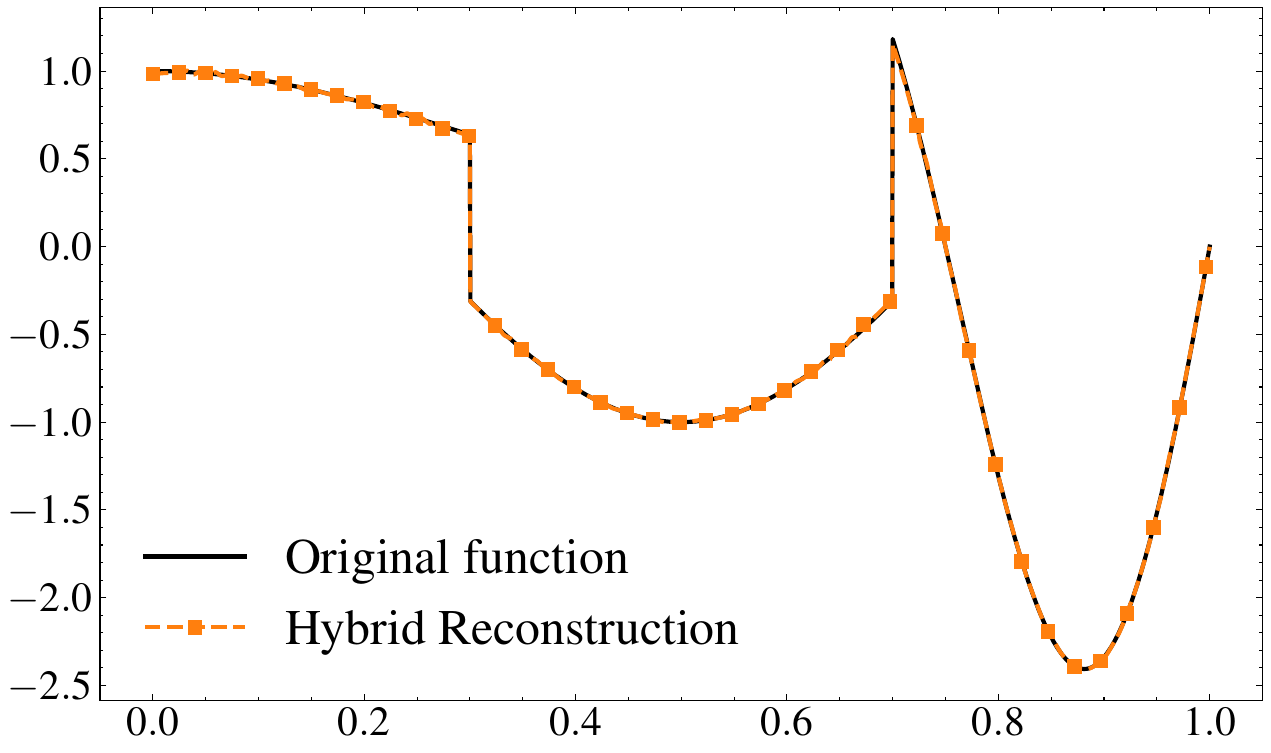}    &
      \includegraphics[width=0.32\textwidth]{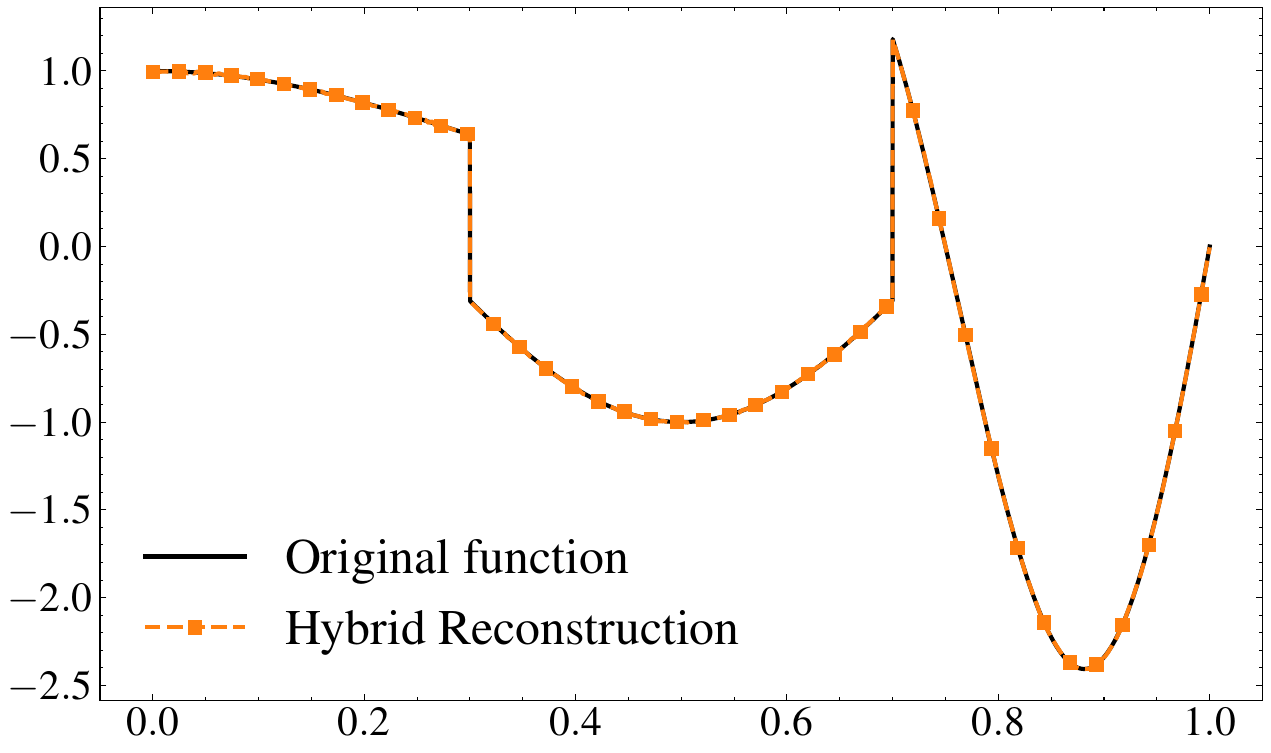}    &
      \includegraphics[width=0.32\textwidth]{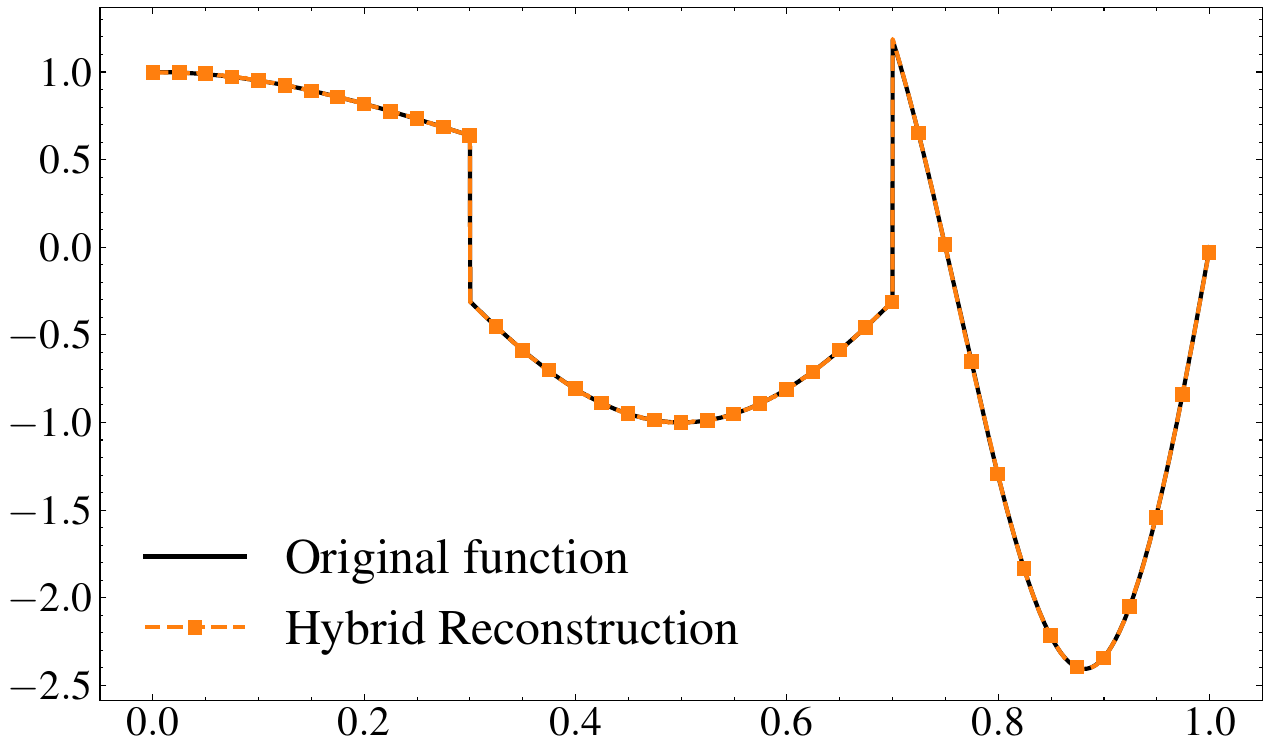}                                      \\
    \end{tabular}
    \caption{Filtered reconstruction \(f_{m}^{\filter}\) (top row) and hybrid reconstruction \(f_{m, M, \delta}^{\hybrid}\) (bottom row) for the jittered sampling with \(m=128, 256, 512\).}
    \label{fig:jittered_reconstruction_multi}
  \end{figure}

  \begin{figure}[!h]
    \centering
    \begin{tabular}{ccc}
      \textbf{m=128}                                                                   & \textbf{m=256} & \textbf{m=512} \\
      \includegraphics[width=0.32\textwidth]{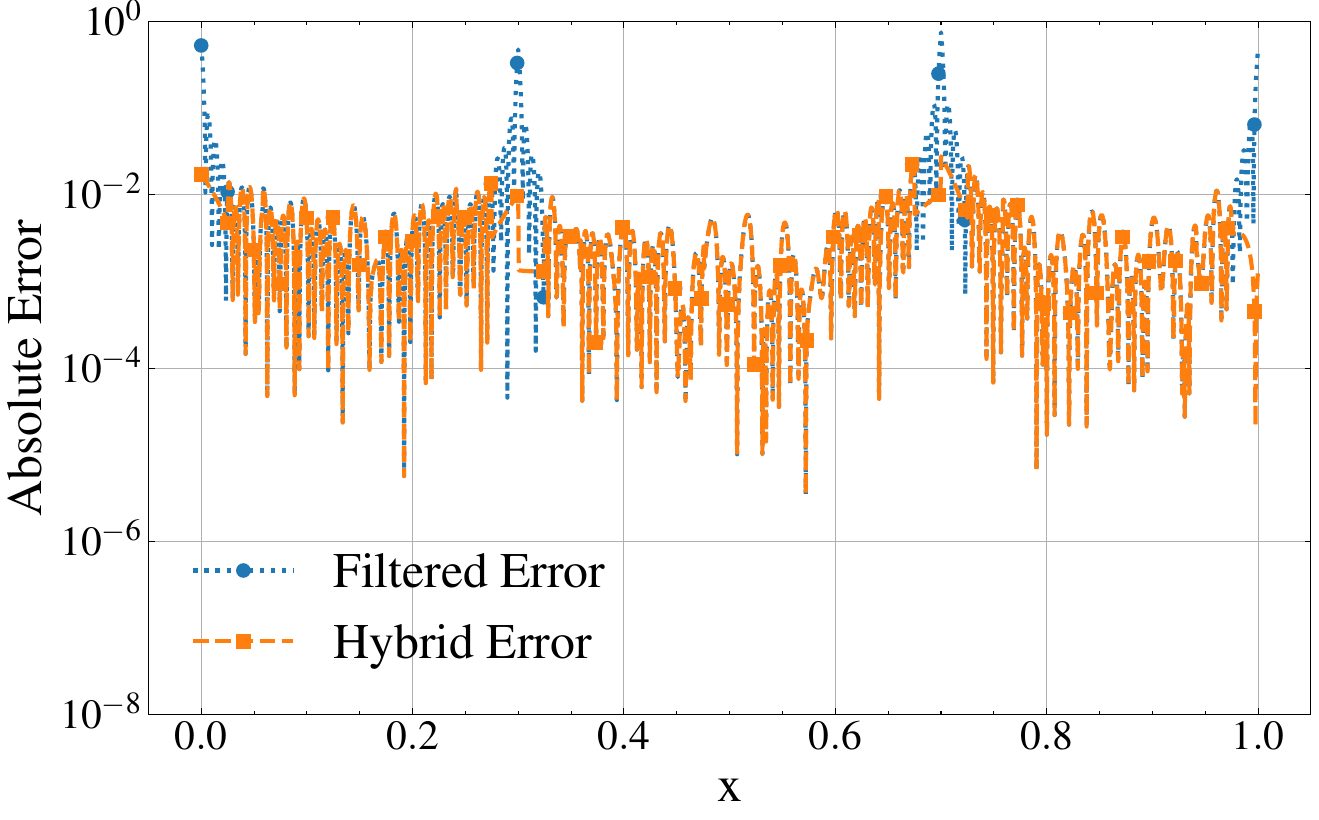} &
      \includegraphics[width=0.32\textwidth]{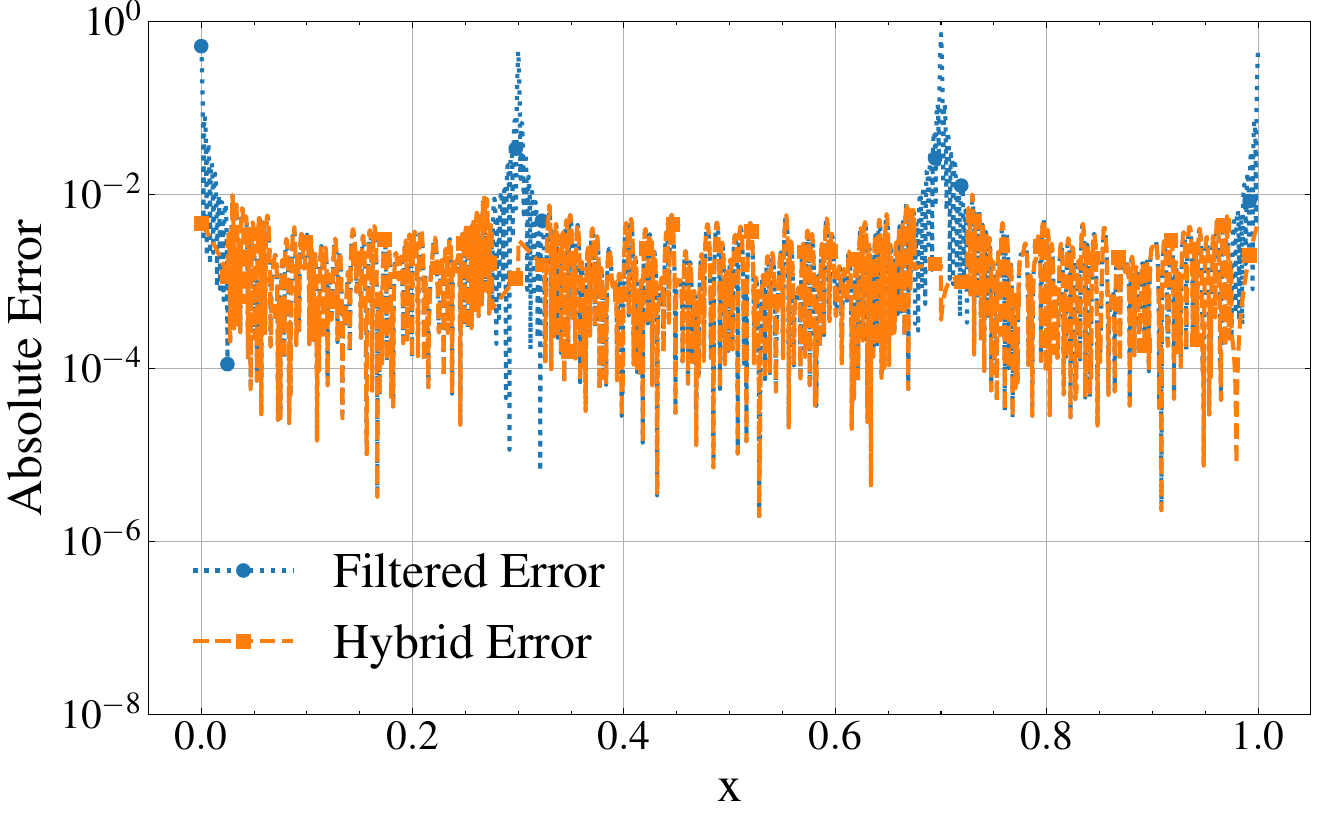} &
      \includegraphics[width=0.32\textwidth]{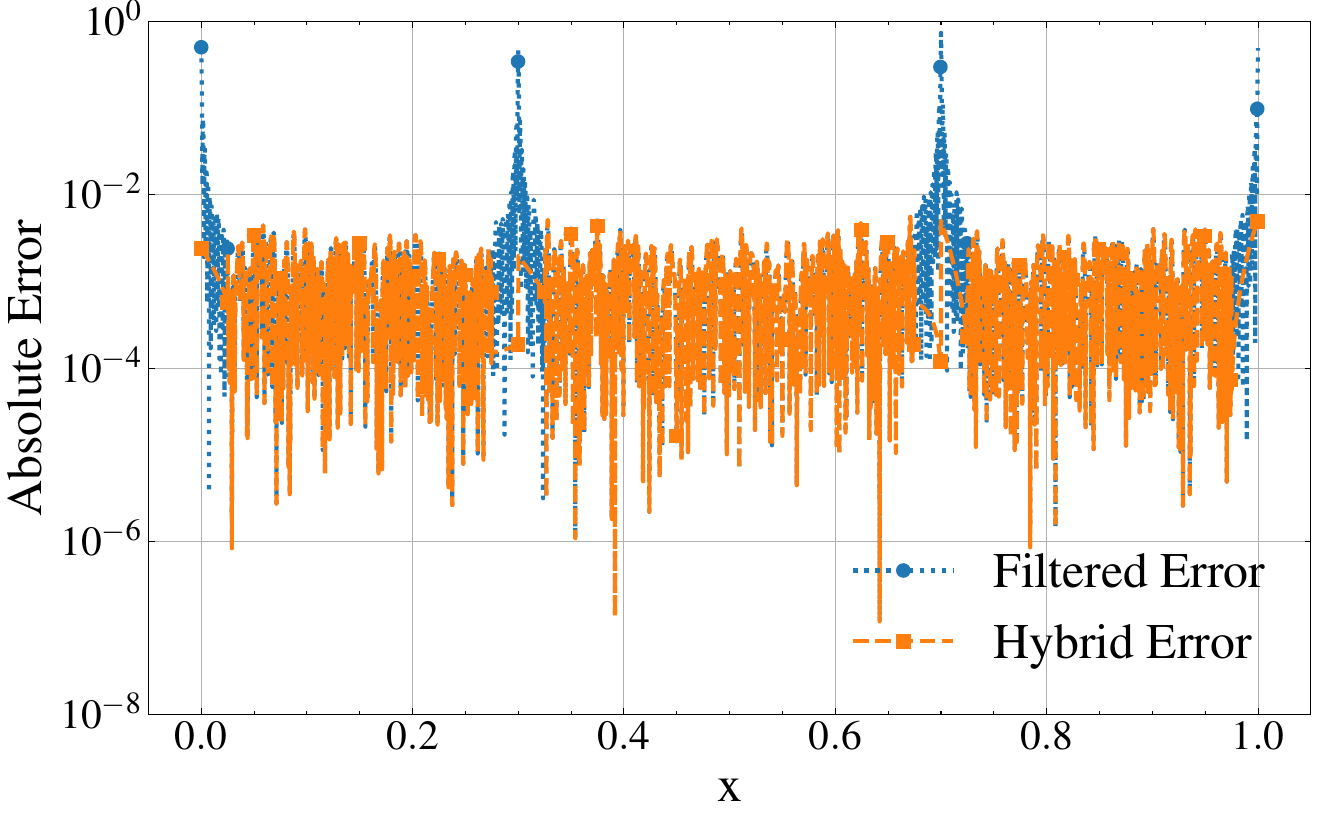}                                   \\
    \end{tabular}
    \caption{Pointwise approximation errors of the filter reconstruction \(f_{m}^{\filter}\)  and the hybrid reconstruction \(f_{m, M, \delta}^{\hybrid}\) for the jittered sampling with \(m=128, 256, 512\).}
    \label{fig:jittered_error_multi}
  \end{figure}

  For log sampling, we display the reconstructed functions in \Cref{fig:log_reconstruction_multi} and the pointwise approximation errors in \Cref{fig:log_error_multi}.
  \begin{figure}[!h]
    \centering
    \begin{tabular}{ccc}
      \textbf{m=128}                                                                        & \textbf{m=256} & \textbf{m=512} \\
      \includegraphics[width=0.32\textwidth]{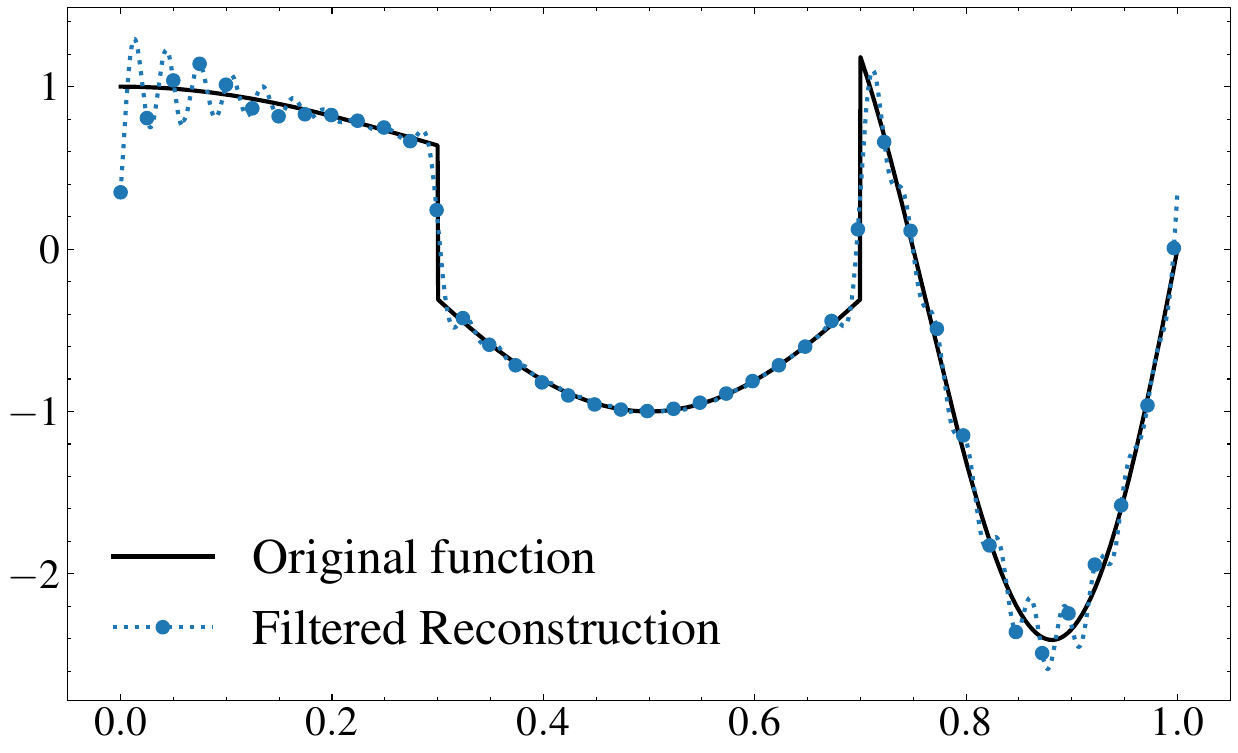} &
      \includegraphics[width=0.32\textwidth]{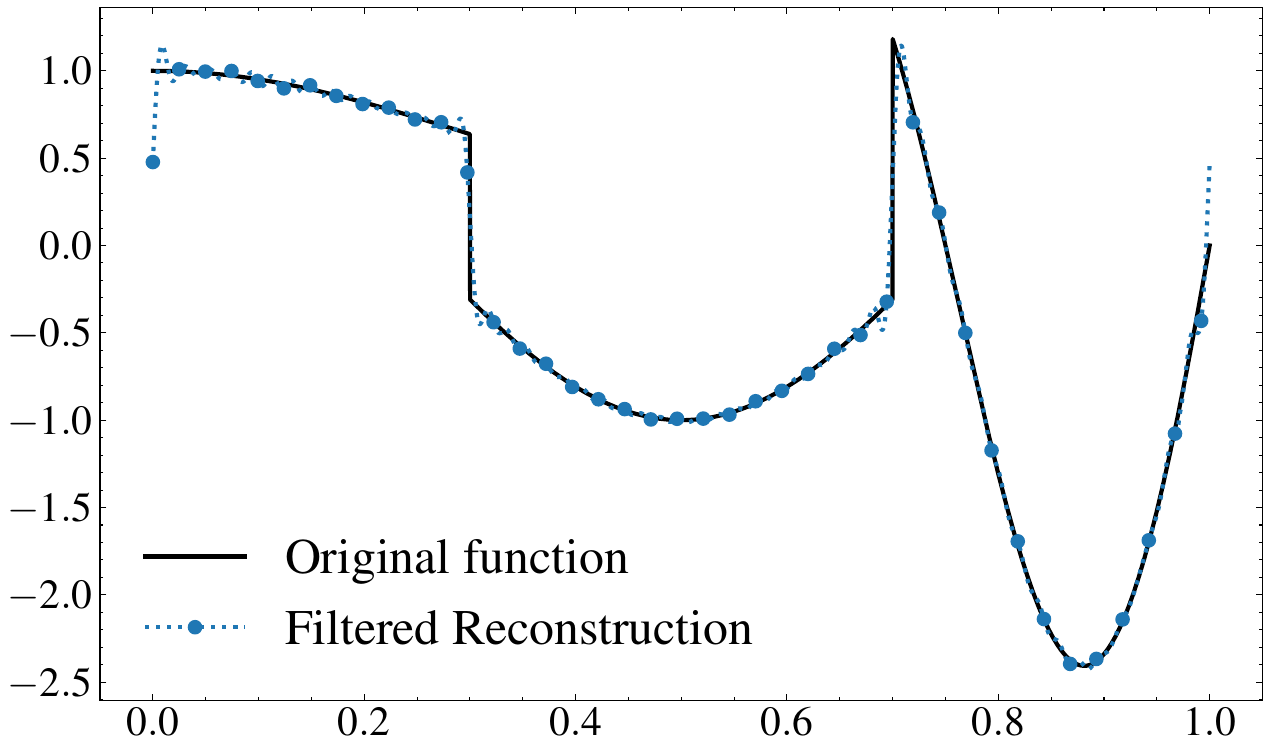} &
      \includegraphics[width=0.32\textwidth]{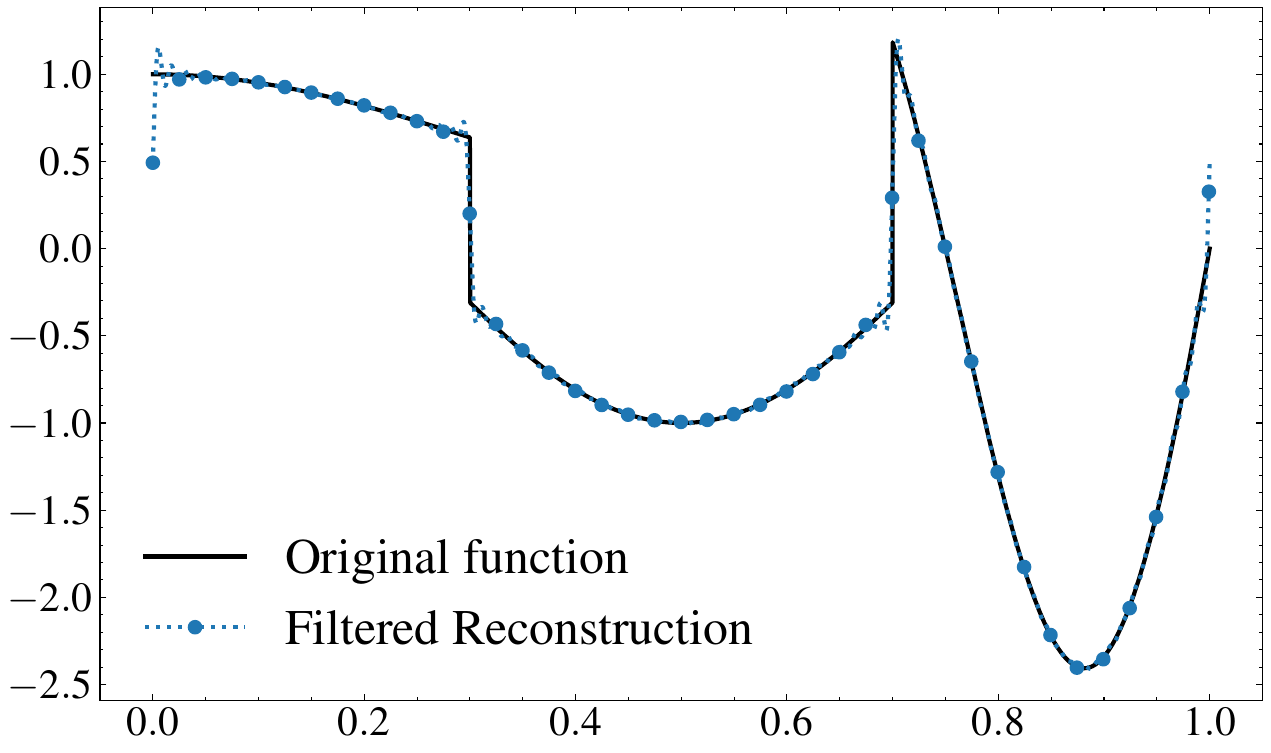}                                   \\
      \includegraphics[width=0.32\textwidth]{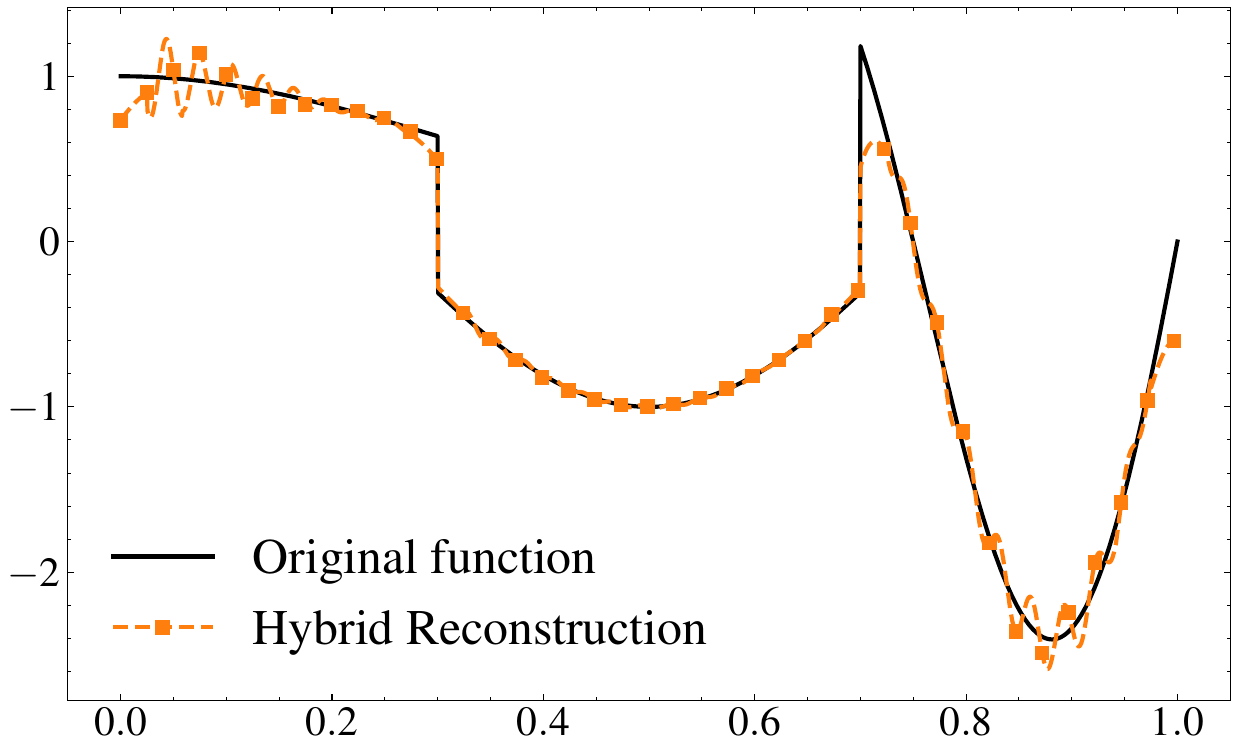}    &
      \includegraphics[width=0.32\textwidth]{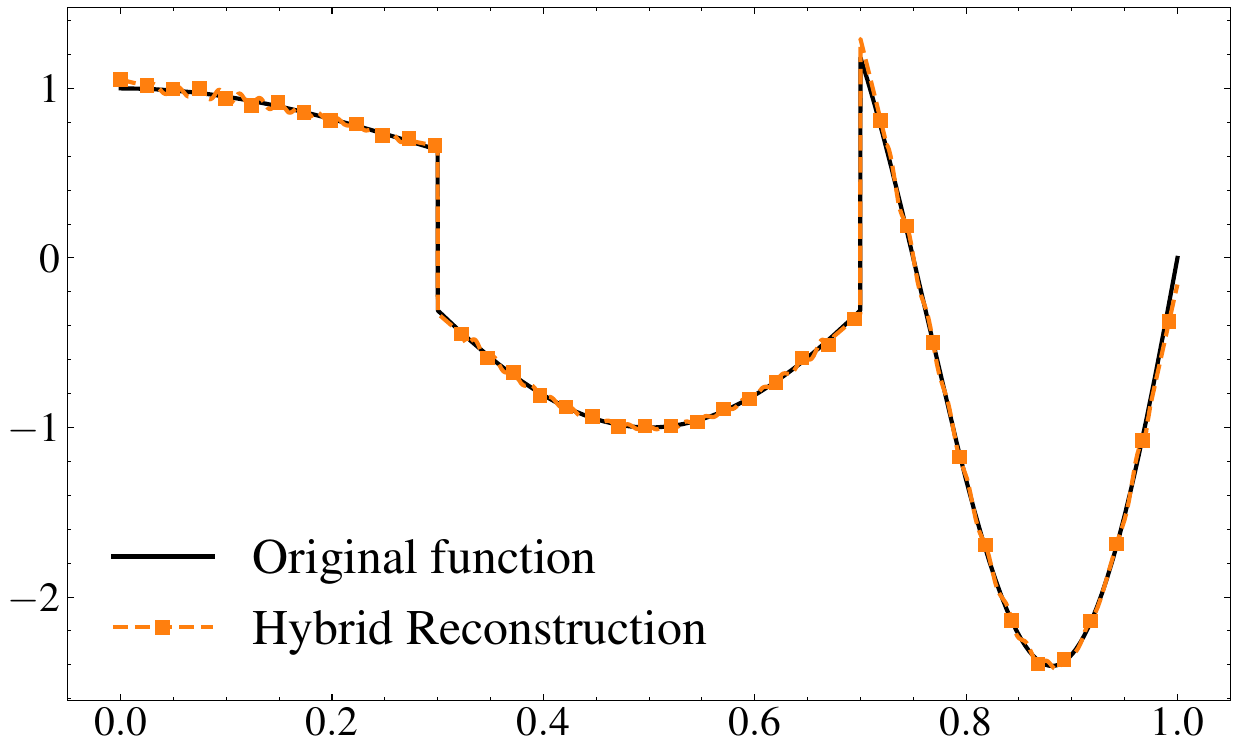}    &
      \includegraphics[width=0.32\textwidth]{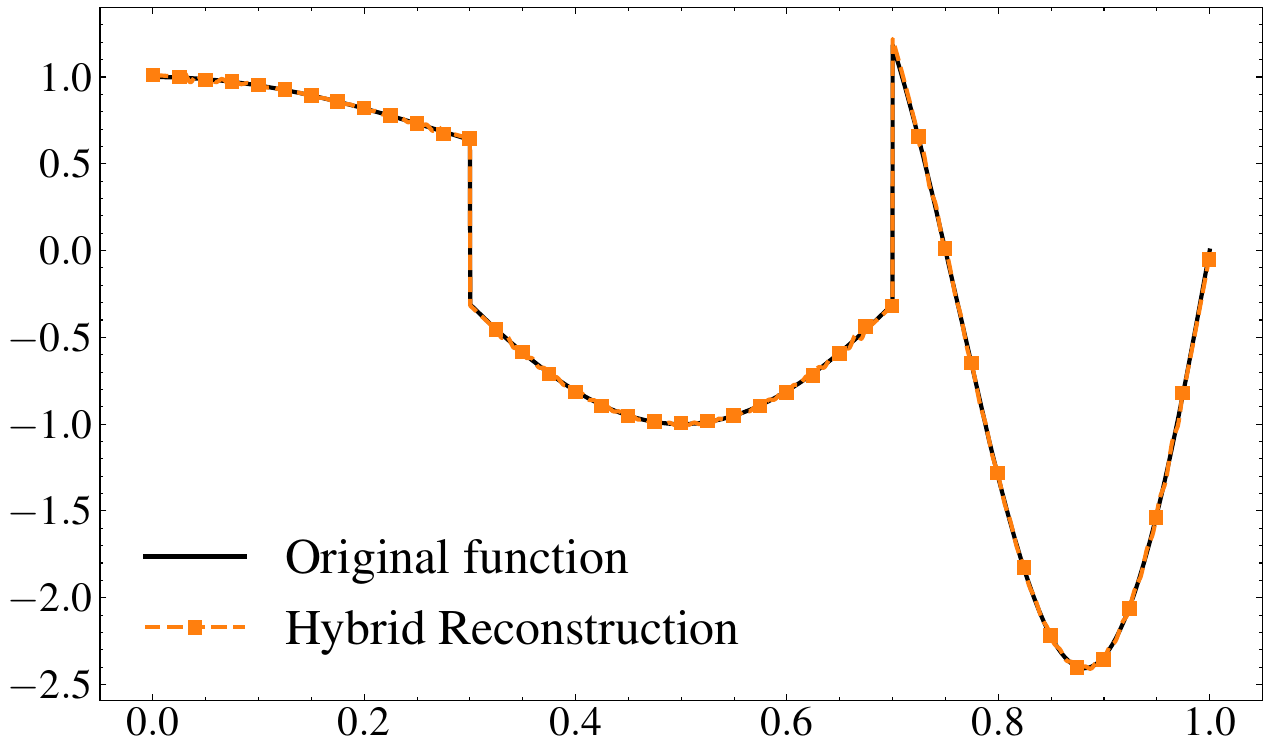}                                      \\
    \end{tabular}
    \caption{Filtered reconstruction \(f_{m}^{\filter}\) (top row) and hybrid reconstruction \(f_{m, M, \delta}^{\hybrid}\) (bottom row) for the log sampling with \(m=128, 256, 512\).}

    \label{fig:log_reconstruction_multi}
  \end{figure}

  \begin{figure}[!h]
    \centering
    \begin{tabular}{ccc}
      \textbf{m=128}                                                                   & \textbf{m=256} & \textbf{m=512} \\
      \includegraphics[width=0.32\textwidth]{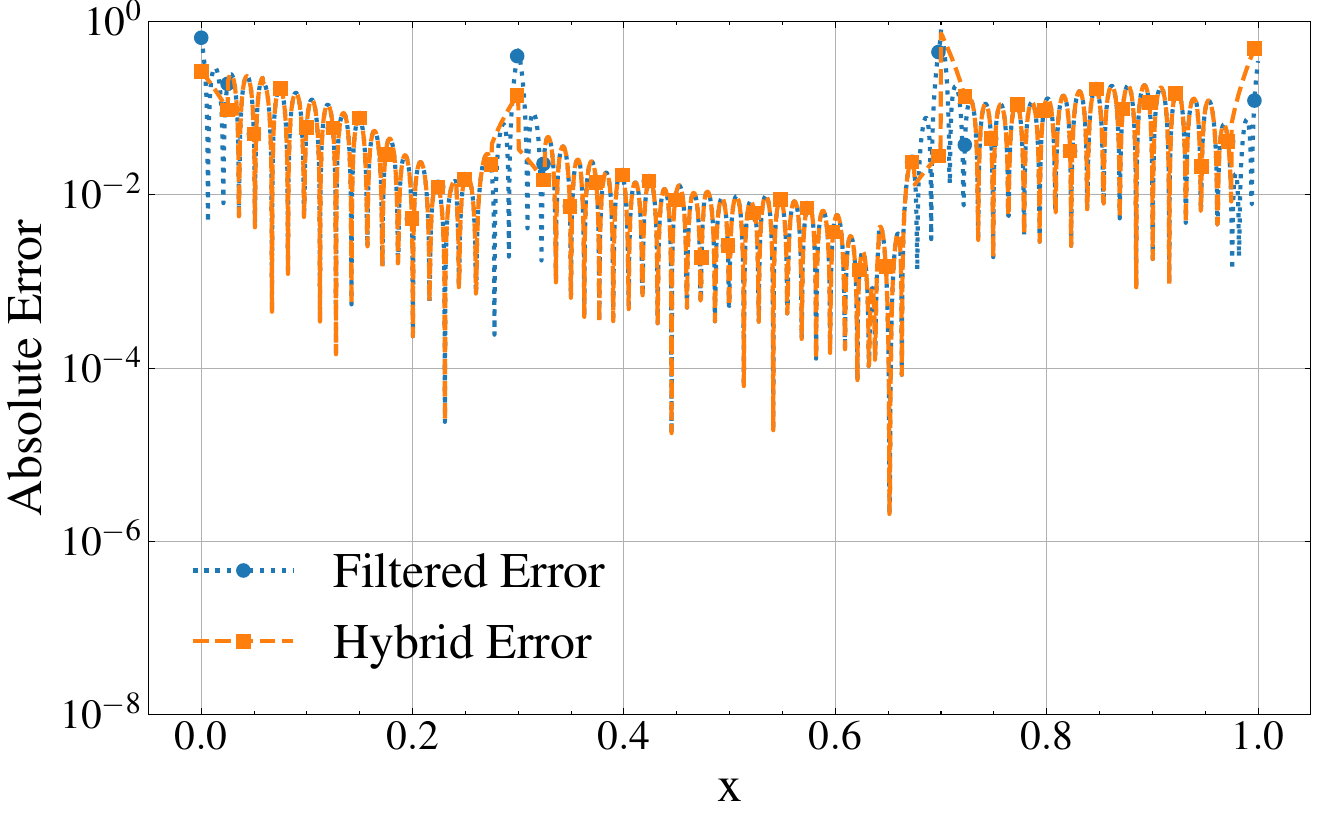} &
      \includegraphics[width=0.32\textwidth]{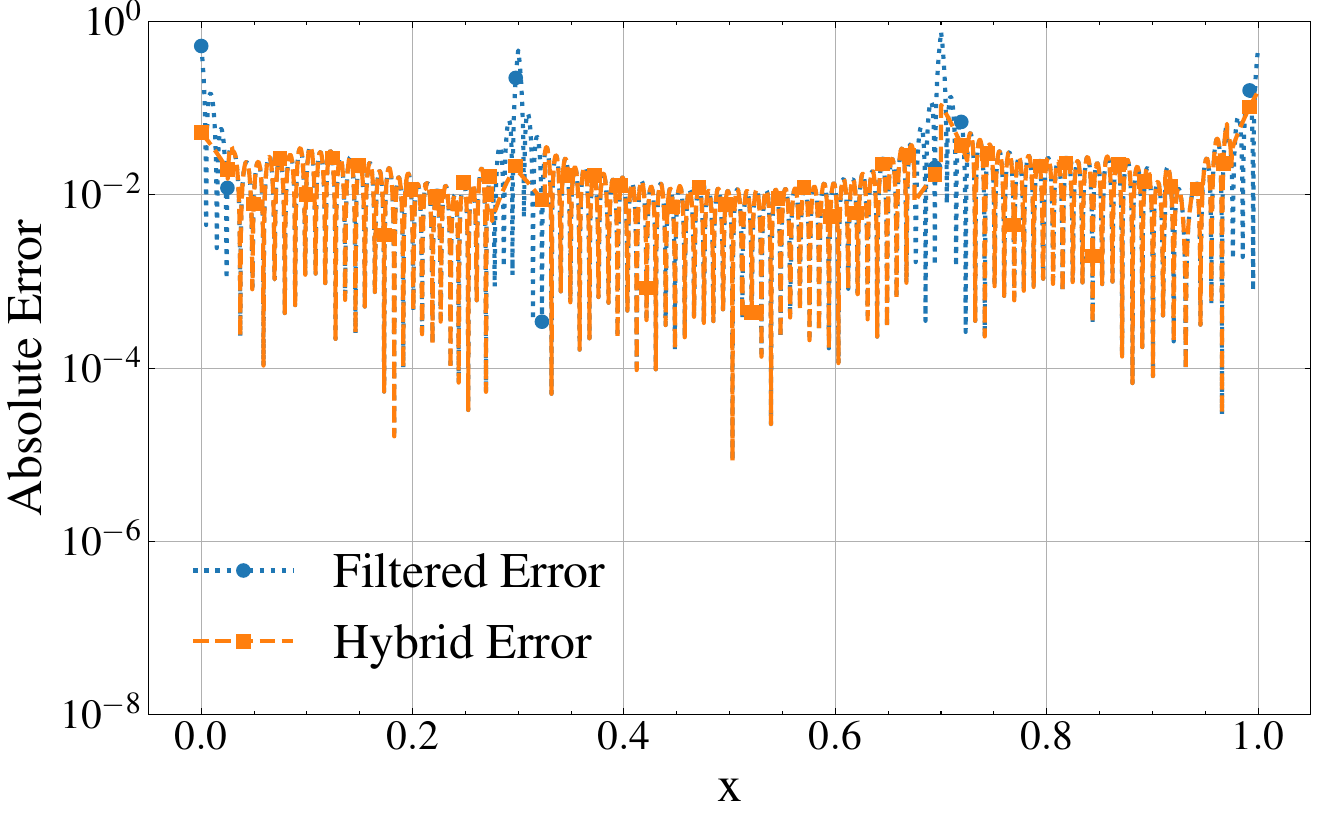} &
      \includegraphics[width=0.32\textwidth]{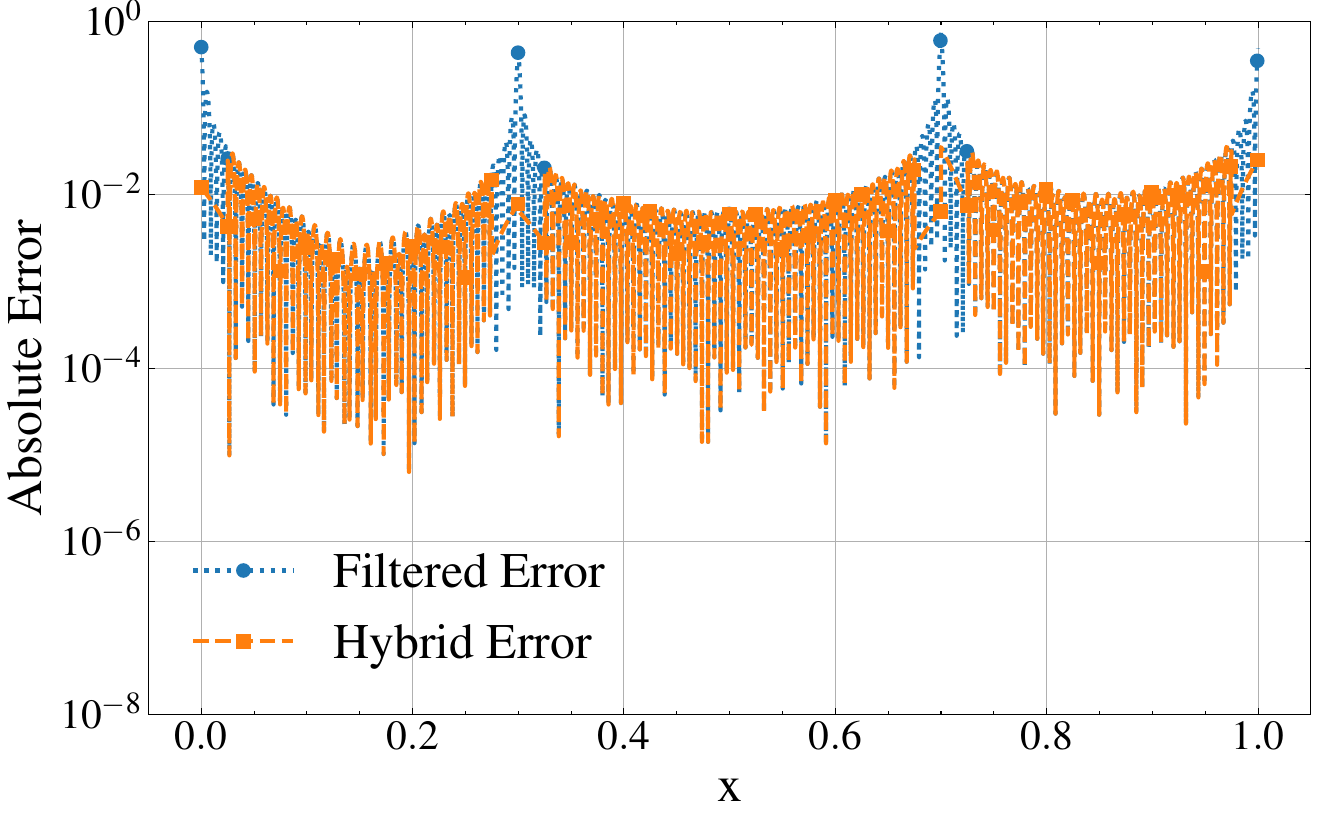}                                   \\
    \end{tabular}
    \caption{Pointwise approximation errors of the filter reconstruction \(f_{m}^{\filter}\)  and the hybrid reconstruction \(f_{m, M, \delta}^{\hybrid}\) for the log sampling with \(m=128, 256, 512\).}
    \label{fig:log_error_multi}

  \end{figure}

\end{example}

We observe in both examples that the filter reconstruction \(f_{m}^{\filter}\) could achieve high-order accuracy away from the jump discontinuities but suffers from low-order accuracy near the jump discontinuities. On the other hand, the proposed hybrid reconstruction \(f_{m, M, \delta}^{\hybrid}\) could achieve high-order accuracy through the whole region \([0,1]\). This demonstrates the effectiveness of the proposed hybrid filter-extrapolation method for reconstructing piece-wise smooth functions from their non-uniform Fourier data.

All the data used to produce the numerical results in this paper are simulated and could be requested from the corresponding author.

\section{Conclusions}\label{sec:conclusion}
In this paper, we have proposed a hybrid filter-extrapolation method for reconstructing a piece-wise smooth function from its non-uniform Fourier data. The proposed method first uses the admissible frame approximation to compute the filter reconstruction of the underlying function, which could achieve high-order accuracy away from the jump discontinuities. It then employs the stable polynomial extrapolation to recover the function values near the jump discontinuities based on the function values of the filter reconstruction away from the jump discontinuities. We have shown that the proposed hybrid method could achieve uniform exponential accuracy through the whole region. Numerical experiments have been presented to demonstrate the performance of the proposed method.

\section*{Declarations}
The project is partially supported by the National Science Foundation DMS-2318781.

\bibliographystyle{siam}
\bibliography{hybrid}

\end{document}